\documentclass[11pt]{elsarticle}

\usepackage{graphics,graphicx}
\usepackage{subfigure}
\usepackage{geometry}

\usepackage{amssymb, amsfonts, amsbsy, amsmath,mathrsfs} 
\usepackage{eucal}

\usepackage{multirow}
\usepackage{lscape}
\usepackage{setspace}

\usepackage{pgfplots}

\newcommand{\Eref}[1]{Equation (\ref{#1})}
\newcommand{\Erefs}[1]{Equations (\ref{#1})}

\newcommand{\fref}[1]{Figure~\ref{#1}}
\newcommand{\frefs}[1]{Figures~\ref{#1}}

\newcommand{\ff}{\mathbf{f}}
\newcommand{\KK}{\mathbf{K}}
\newcommand{\bm}{\mathbf{M}}

\newcommand{\bveps}{\boldsymbol{\varepsilon}}
\newcommand{\bvsig}{\boldsymbol{\sigma}}

\usepackage{setspace}
\begin{document}

\begin{frontmatter}
%\title{Static bending and free vibration analysis of sandwich plates with nanotube reinforced face sheets using a higher order accurate theory}

%\title{A comparison of various plate theories for sandwich plates with homogeneous core and nanotube reinforced facesheets}
\title{Application of higher-order structural theory to  bending and free vibration analysis of sandwich plates with CNT reinforced composite facesheets}

\author[a]{Sundararajan Natarajan\fnref{label2}\corref{cor1}}
\author[b]{Mohamed Haboussi}
%\author[a]{Fan Jiang}
\author[c]{Manickam Ganapathi}

\address[a]{School of Civil and Environmental Engineering, The University of New South Wales, Sydney, NSW 2052, Australia.}
\address[b]{Universit\'e Paris 13-CNRS, LSPM, UPR 3407, Villetaneuse, F-93430, France.}
%\address[c]{Head, Stress \& DTA, IES-Aerospace, Mahindra Satyam Computers Services Ltd., Bangalore, India}
\address[c]{Tech Mahindra Ltd., Electronic City, Bangalore-560100, India.}

\fntext[3]{School of Civil and Environmental Engineering, The University of New South Wales, Sydney, NSW 2052, Australia. Tel: +61 293855030; E-mail: s.natarajan@unsw.edu.au; snatarajan@cardiffalumni.org.uk}

\begin{abstract}
In this paper, the bending and free flexural vibration behaviour of sandwich plates with carbon nanotube (CNT) reinforced facesheets are investigated using QUAD-8 shear flexible element developed based on higher-order structural theory. This theory accounts for the realistic variation of the displacements through the thickness, and the possible discontinuity in slope at the interface, and the thickness stretch affecting the transverse deflection. The in-plane and rotary inertia terms are considered in the formulation. The governing equations obtained using Lagrange's equation of motions are solved for static and dynamic analyses considering a sandwich plate with homogeneous core and CNT reinforced face sheets. The accuracy of the present formulation is tested considering the problems for which solutions are available. A detailed numerical study is carried out based on various higher-order models deduced from the present theory to examine the influence of the volume fraction of the CNT, core-to-face sheet thickness and the plate thickness ratio on the global/local response of different sandwich plates.
\end{abstract}

\begin{keyword} 
Carbon nanotube reinforcement \sep sandwich plate higher-order theory \sep mechanical loading \sep thermal loading \sep free vibration \sep shear flexible element.
\end{keyword}

\end{frontmatter}

\section{Introduction}
Engineered materials, mostly inspired from nature have been replacing the conventional materials due to their high strength-to- and stiffness-to-weight ratios. Among the various structural constructions, the sandwich type of structures are more attractive due to its outstanding bending rigidity, low specific weight, excellent vibration characteristics and good fatigue properties. These sandwich constructions can be considered for the requirement of light weight and high bending stiffness in design by appropriate choice of materials. The response of such structures depends on the bonding characteristics. A typical sandwich structure may consist of a homogeneous core with face sheets. To improve the characteristics, the face sheets can be laminated composites~\cite{whitney1972}, functionally graded materials~\cite{Zenkour2005} or polymer matrix with reinforcements~\cite{ugalesingh2013}. Recent introduction of carbon nanotubes (CNTs) as reinforcement has attracted researchers to investigate the responses of such structures~\cite{tjong2009}. Experimental investigations show that the CNTs have extraordinary mechanical properties than those of carbon fibers~\cite{jiazhao2011}. The CNTs are seen as promising reinforcement for sandwich face sheets. Thostenson and Chou~\cite{thostensonchou2002} showed that the addition of nanotubes increases the tensile modulus, yield strength and ultimate strengths of the polymer films. Their study has also showed that the polymer films with aligned nanotubes as reinforcements yield superior strength to randomly oriented nanotubes. The properties of the polymer films can also be optimized by varying the distribution of CNTs through the thickness of the film. Moreover, Formica \textit{et al.,}~\cite{formicalacarbonara2010} highlighted that the CNT reinforced plates can be be tailored to respond to an external excitation. This has generated interest among researchers. For predicting the realistic behaviour of sandwich structures, more accurate analytical/numerical models based on the three-dimensional models may be computationally involved and expensive. Hence, among the researchers, there is a growing appreciation of the importance of applying two-dimensional theories with new kinematics for the evolution of the accurate structural analysis. Various structural theories proposed for the laminated structures have been examined and some of the important contributions pertaining to the sandwich laminated plates are discussed here. 

Based on the first-order shear deformation theory, Zhu \textit{et al.,}~\cite{zhulei2012}  studied the static and free vibration of CNT reinforced plates. They considered polymer matrix with CNT reinforcement, neglecting the temperature effects. It was predicted that the CNT volume fraction has greater influence on the fundamental frequency and the maximum center deflection. Wang and Shen~\cite{wangshen2011} studied the large amplitude vibration of nanocomposite plates resting on elastic foundation using a perturbation technique. The governing equations were based on higher-order shear deformation theory and the composite plates were reinforced with carbon nanotubes. Their study brought out that while the linear frequencies decrease with the addition of CNTs, the nonlinear to linear frequency ratio increased, especially when increasing the temperature or by decreasing the foundation stiffness. Arani \textit{et al.,}~\cite{aranimaghamikia2011}, Liew \textit{et al.,}~\cite{liewlei2014} and Lei \textit{et al.,}~\cite{leiliew2013} studied the buckling and post-buckling characteristics of CNT reinforced plates using the finite element method and meshless methods, respectively. It was shown that the reinforcement with CNT increasing the load carrying capacity of the plate. Aragh \textit{et al.,}~\cite{araghbarati2012} used generalized differential quadrature and obtained a semi-analytical solution for 3D vibration of cylindrical panels. It was shown that graded CNTs with symmetric distribution through the thickness have high capabilities to alter the natural frequencies when compared to the uniformly distributed or asymmetrically distributed.

It is observed from these investigations that first- order shear deformation theory has been widely employed by many researchers whereas higher-order model considering variation in in-plane displacements has been used by few authors for the analysis of CNT reinforced plates. However, the available literature pertaining to sandwich structures with CNT reinforced face sheets is rather limited compared to plates. Also, to the author’s knowledge, the theories accounting the variation of in-plane displacement through the thickness, and the possible discontinuity in slope at the interface, and the thickness stretch affecting transverse deflection is not exploited while investigating the structural behavior of CNT reinforced sandwich structures. A layerwise theory is the possible candidature for this purpose, but it may be computationally expensive as the number of unknowns to be solved increases with increase in the number of mathematical or physical layers. Ali \textit{et al.,}~\cite{Ali1999} and Ganapathi and Makhecha~\cite{Ganapathi2001} have used an alternative higher-order plate theory based on global approach, for multi-layered laminated composites by incorporating the realistic through the thickness approximations of the in-plane and transverse displacements by adding a zig-zag function and higher-order terms, respectively. This formulation has proved to give very accurate results for the composite laminates. Such a model for the current problems is considered while comparing with the other approaches available in the literature

%Such model for the current problems may be worthwhile considering as a candidature while comparing with the other approaches available in the literature.

\paragraph{Approach} In this paper, a $\mathcal{C}^o$ 8-noded quadrilateral plate element with 13 degrees of freedom per node ~\cite{Ganapathi2001,Makhecha2001,natarajanmanickam2012} based on higher order theory~\cite{Ali1999} is employed to study the static deflection and free vibration analysis of thick/thin sandwich carbon nanotube reinforced functionally graded material plates. The efficacy of the present formulation, for the static analyses subjected to mechanical and thermal loads, and the free vibration study is illustrated through the numerical studies by employing various structural models deduced from the present higher-order accurate theory.

\paragraph{Outline} The paper is organized as follows. The computation of the effective properties of carbon nanotube reinforced composites are discussed in the next section. Section \S \ref{highertheory} presents the higher order accurate theory to describe the plate kinematics and Section \S \ref{eledes} describes the 8-noded quadrilateral plate element employed in this study. The numerical results for the static deflection and the free vibration of thick/thin sandwich carbon nanotube reinforced functionally graded plates are given in Section \S \ref{numexamples}, followed by concluding remarks in the last section.

\section{Theoretical Formulation}
\label{formulate}
\begin{figure}[htpb]
\centering
\scalebox{0.8}{\input{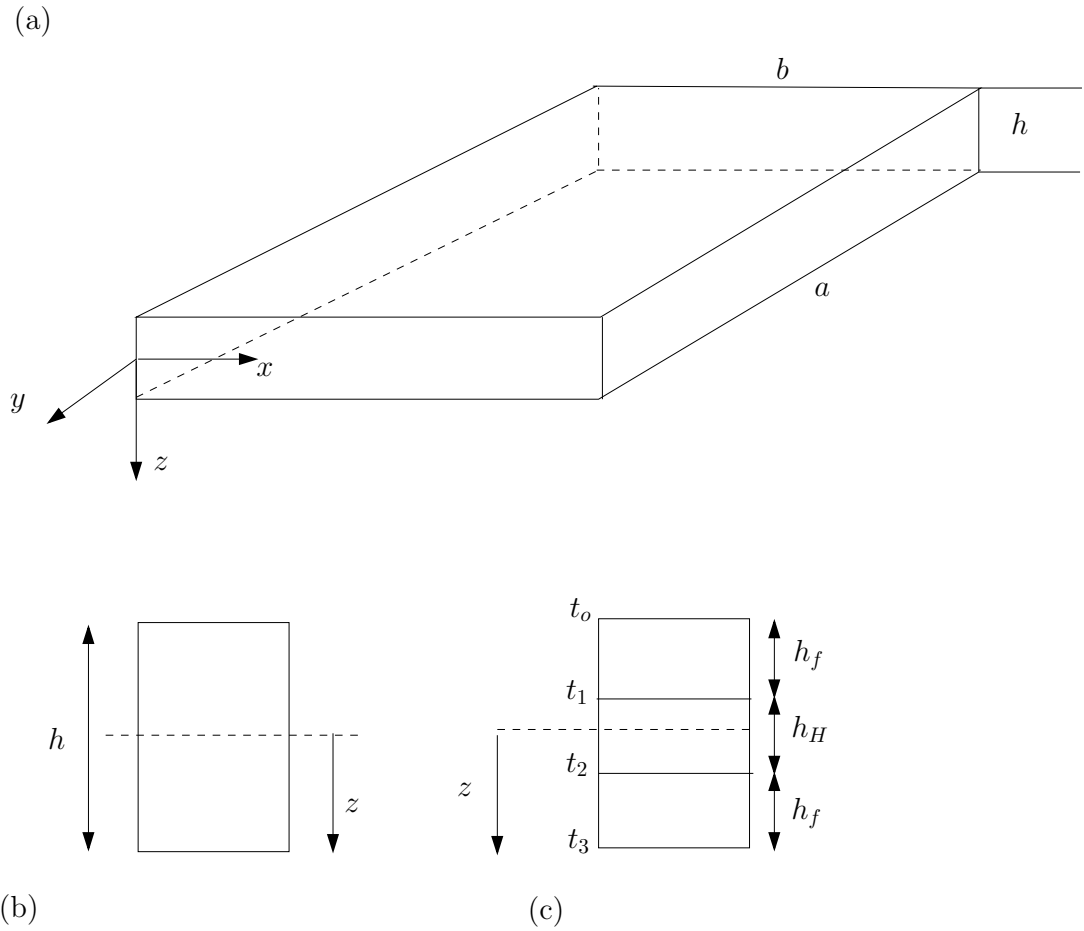}}
\caption{Coordinate system of a rectangular plate: (a) $x$ and $y$ along the in-plane directions and $z$ along the thickness direction: (b) single layer plate cross section and (c) sandwich plate made up of three layers, where $h_H$ is the core thickness and $h_f$ is the facesheet thickness. The CNTs are functionally graded in the facesheets.}
\label{fig:platedescript}
\end{figure}

Consider a rectangular sandwich plate with co-ordinates $x$ and $y$ along the in-plane directons and $z$ along the thickness direction as shown in \fref{fig:platedescript}. The core-to-facesheet thickness ratio is $h_H/h_f$, where, $h_H$ is the core thickness and $h_f$ is the facesheet thickness. The length, the width and the total thickness of the plate are $a, b$ and $h$ (see \fref{fig:platedescript}). A Cartesian coordinate system is assumed and the origin is located at the corner of the plate on the middle plane. We assume that the CNT reinforced layer is made from a mixture of single walled CNT, graded distribution in the thickness direction and the matrix is assumed to be isotropic. The effective properties of such reinforced structures can be computed by Mori-Tanaka scheme~\cite{araghbarati2012} or by the rule of mixtures. In this study, we employ simple rule of mixtures with correction factors to estimate the effective material properties of CNT reinforced matrix. The effective material properties of the CNT reinforced matrix are given by~\cite{wangshen2012}:
\begin{align}
E_{11} &= \eta_1 V_{\rm CN} E_{11}^{\rm CN} + V_m E_m \nonumber \\
\frac{\eta_2}{E_{22}} &= \frac{ V_{\rm CN}}{E_{22}^{\rm CN}} + \frac{V_m}{E_m} \nonumber \\
\frac{\eta_3}{G_{12}} &= \frac{V_{\rm CN}}{G_{12}^{\rm CN}} + \frac{V_m}{G_m} \nonumber \\
\nu_{12}^{\rm CN} &= \nu_{12}^{\rm CN} V_{\rm CN}^\ast + \nu_m V_m \nonumber \\
\rho &= \rho_{\rm CN} V_{\rm CN} + \rho_m V_m
\label{eqn:effecprop}
\end{align}
where $E_{11}^{\rm CN}, E_{22}^{\rm CN}$ and $G_{12}^{\rm CN}$ are the Young's modulii and the shear modulus of the CNT, respectively and $V_{\rm CN}$ and $V_m$ are the volume fraction of the CNT and the matrix, respectively. The volume fractions are related by: $V_{\rm CN} + V_m=1$. The CNT efficiency parameters $(\eta_1,\eta_2,\eta_3)$ are introduced to account for the inconsistency in the load transfer between the CNT and the matrix. The values of the efficiency parameters are obtained by matching the elastic modulus of CNT reinforced polymer matrix observed from the MD simulation results with the numerical results obtained from the rule of mixtures. 

\paragraph{Three layer sandwich structure} The three layer sandwich structure is made up of a homogeneous core with core thickness $h_H$ and two facesheets with thickness $h_f$ (see \fref{fig:platedescript}). The facesheets are assumed to be reinforced with CNTs. We assume the volume fraction $V_{\rm CN}$ for the top facesheet as
\begin{equation}
V_{\rm CN} = 2 \left( \frac{t_1 - z}{t_1-t_0} \right) V_{\rm CN}^\ast
\label{eqn:topsh}
\end{equation}
and for the bottom facesheet as follows
\begin{equation}
V_{\rm CN} = 2 \left( \frac{z-t_2}{t_3-t_2} \right) V_{\rm CN}^\ast
\label{eqn:botsh}
\end{equation}
where,
\begin{equation*}
V_{\rm CN}^\ast = \frac{w_{\rm CN}}{w_{\rm CN} + \left( \frac{\rho_{\rm CN}}{\rho_m} \right) \left[1 - w_{\rm CN} \right]}
\end{equation*}
where $w_{\rm CN}$ is the mass fraction of the nanotube, $\rho_{\rm CN}$ and $\rho_m$ are the mass densities of carbon nanotube and the matrix, respectively. The thermal expansion coefficient in the longitudinal and transverse directions can be expressed as~\cite{wangshen2012}:
\begin{align}
\alpha_{11} &= \alpha_{11}^{\rm CN} V_{\rm CN} + \alpha_m V_m \nonumber \\
\alpha_{22} &= (1+\nu_{12}^{\rm CN})V_{\rm CN}\alpha_{22}^{\rm CN} + (1+\nu_m)V_m\alpha_m - \nu_{12}\alpha_{11}
\end{align}
where $\alpha_{11}^{\rm CN}, \alpha_{22}^{\rm CN}$ and $\alpha_m$ are the thermal expansion coefficients for the CNT and the matrix, respectively and $\nu_{12}^{\rm CN}$ and $\nu_m$ are the Poisson's ratio.

%--------- Higher order accurate theory
\section{Higher order accurate theory} \label{highertheory}
The sandwich plate is assumed to be made up of three discrete layers with a homogeneous core. The in-plane displacements $u_k$ and $v^k$, and the transverse displacement $w^k$ for the $k^{th}$ layer, are assumed as~\cite{Ali1999,Makhecha2001}:
\begin{eqnarray}
u^k(x,y,z,t) = u_o(x,y,z,t) + z \theta_x(x,y,t) + z^2 \beta_x(x,y,t) + z^3 \phi_x(x,y,t) + S^k \psi_x(x,y,t)  \nonumber \\
v^k(x,y,z,t) = v_o(x,y,z,t) + z \theta_y(x,y,t) + z^2 \beta_y(x,y,t) + z^3 \phi_y(x,y,t) + S^k \psi_y(x,y,t)  \nonumber \\
w^k(x,y,z,t) = w_o(x,y,t) + z w_1(x,y,t) + z^2 \Gamma(x,y,t)
\label{eqn:dispField}
\end{eqnarray}
The terms with even powers of $z$ in the in-plane displacements and odd powers of $z$ occurring in the expansion for $w^k$ correspond to the stretching problem. However, the terms with odd powers of $z$ in the in-plane displacements and the even ones in the expression for $w^k$ represent the flexure problem. $u_o, v_o$ and $w_o$ are the displacements of a generic point on the reference surface; $\theta_x$ and $\theta_y$ are the rotations of the normal to the reference surface about the $y$ and $x$ axes, respectively; $w_1,\beta_x,\beta_y,\Gamma,\phi_x$ and $\phi_y$ are the higher-order terms in the Taylor`s series expansions, defined at the reference surface. $\psi_x$ and $\psi_y$ are generalized variables associated with the zigzag function, $S^k$. The zigzag function, $S^k$, as given in~\cite{Murukami1986,Ganapathi2001,Roderigues2011}, is defined by
\begin{equation}
S^k = 2(-1)^k \frac{z_k}{h_k}
\end{equation}
where $z_k$ is the local transverse coordinate with the origin at the center of the $k^{th}$ layer and $h_k$ is the corresponding layer thickness. Thus, the zigzag function is piecewise linear with values of -1 and 1 alternatively at different interfaces. The `zig-zag' function, as defined above, takes care of the inclusion of the slope discontinuities of $u$ and $v$ at the interfaces of the sandwich plate as observed in the exact three-dimensional elasticity solutions of thick sandwich functionally graded materials. The main advantage of using such formulation is that such a function is more economical than a discrete layer approach~\cite{Nosier1993,Ferreira2005}.

The strains in terms of mid-plane deformation, rotations of normal and higher order terms associated with displacements are:
\begin{equation}
\bveps = \left \{ \begin{array}{c} \bveps_{\rm bm} \\ \bveps_s \end{array} \right\}
\label{eqn:strain}
\end{equation}
The vector  $\bveps_{\rm bm}$ includes the bending and membrane terms of the strain components and vector $\bveps_s$ contains the transverse shear strain terms. These strain vectors are defined as:
\begin{eqnarray}
\bveps_{\rm bm} &=& \left\{ \begin{array}{c} \varepsilon_{xx} \\ \varepsilon_{yy} \\ \varepsilon_{zz} \\ \gamma_{xy} \end{array} \right\} + \left\{ \begin{array}{c} u_{,x} \\ v_{,y} \\ w_{,z} \\ u_{,y} + v_{,x} \end{array} \right\} \nonumber \\
&=& \bveps_0 + z \bveps_1 + z^2 \bveps_2 + z^3 \bveps_3 + S^k \bveps_4 
\end{eqnarray}
\begin{eqnarray}
\bveps_s &=& \left\{ \begin{array}{c} \gamma_{xz} \\ \gamma_{yz} \end{array} \right\} = \left\{ \begin{array}{c} u_{,z} + w_{,x} \\ v_{,z} + w_{,y} \end{array} \right\} \nonumber \\
&=& \gamma_o + z \gamma_1 + z^2 \gamma_2 + S^k_{,z} \gamma_3
\end{eqnarray}
where
\begin{eqnarray}
\bveps_o = \left\{ \begin{array}{c} u_{o,x} \\ v_{o,y} \\ w_1 \\ u_{o,y} + v_{o,x} \end{array} \right\},\hspace{1cm} \bveps_1 = \left\{ \begin{array}{c} \theta_{x,x} \\ \theta_{y,y} \\ 2\Gamma \\ \theta_{x,y} + \theta_{y,x} \end{array} \right\}, \nonumber \\
\bveps_2 = \left\{ \begin{array}{c} \beta_{x,x} \\ \beta_{y,y} \\ 0 \\ \beta_{x,y} + \beta_{y,x} \end{array} \right\}, \hspace{1cm} \bveps_3 = \left\{ \begin{array}{c} \phi_{x,x} \\ \phi_{y,y} \\ 0 \\ \phi_{x,y} + \phi_{y,x} \end{array} \right\}, \nonumber \\
\bveps_4 = \left\{ \begin{array}{c} \psi_{x,x} \\ \psi_{y,y} \\ 0 \\ \psi_{x,y} + \psi_{y,x} \end{array} \right\}.
\end{eqnarray}
and,
\begin{eqnarray}
\gamma_o = \left\{ \begin{array}{c} \theta_x + w_{o,x} \\ \theta_y + w_{o,y} \end{array} \right\}, \hspace{1cm} \gamma_1 = \left\{ \begin{array}{c} 2\beta_x + w_{1,x} \\ 2\beta_y + w_{1,y} \end{array} \right\}, \nonumber \\
\gamma_2 = \left\{ \begin{array}{c} 3\phi_x + \Gamma_{,x} \\ 3\phi_y + \Gamma_{,y} \end{array} \right\}, \hspace{1cm} \gamma_3 = \left\{ \begin{array}{c} \psi_x S_{,z}^k \\ \psi_y S_{,z}^k \end{array} \right\} .
\end{eqnarray}
The subscript comma denotes partial derivatives with respect to the spatial coordinate succeeding it. The constitutive relations for an arbitrary layer $k$ can be expressed as:
\begin{eqnarray}
\bvsig &=& \left\{ \begin{array}{cccccc} \sigma_{xx} & \sigma_{yy} & \sigma_{zz} & \sigma_{xz} & \sigma_{yz} \end{array} \right\}^{\rm T} \nonumber \\
&=& \bf{Q}^k \left\{ \begin{array}{cc} \bveps_{\rm bm} & \bveps_s \end{array} \right\}^{\rm T}
\end{eqnarray}
where $\bf{Q}_k$ is the stiffness coefficient matrix defined as
\begin{eqnarray}
Q_{11}^k  = {E_{11} \over 1-\nu_{12}\nu_{21}}; \hspace{0.5cm} Q_{22}^k  = {E_{22} \over 1-\nu_{12}\nu_{21}}; \hspace{0.5cm} Q_{12}^k = {\nu_{21} E_{11} \over 1-\nu_{12}\nu_{21}};  \nonumber \\
Q_{44}^k = G_{23}; \hspace{0.25cm} Q_{55}^k = G_{13}; \hspace{0.25cm} Q_{66}^k = G_{12};
Q_{16}^k = Q_{26}^k = 0.
\label{eqn:stiffcoeff}
\end{eqnarray}
For the homogeneous core, the shear modulus $G$ is related to the Youngs's modulus by: $E = 2G(1+\nu)$. The governing equations of motion are obtained by applying Lagrangian equations of motion given by
\begin{equation}
\frac{d}{dt} \left[ \frac{\partial (T-U)}{\partial \dot{\delta_i}} \right] - \left[ \frac{\partial (T-U)}{\partial \delta_i} \right] = 0, \hspace{1cm} i = 1,2,\cdots,n
\label{eqn:lagrange}
\end{equation}
where $\delta_i$ is the vector of degrees of freedom and $T$ is the kinetic energy of the plate given by:
\begin{equation}
T(\delta) = \frac{1}{2} \iint \left[ \sum_{k=1}^n \int\limits_{h_k}^{h_{k+1}} \rho_k \left\{ \begin{array}{ccc} \dot{u}_k & \dot{v}_k & \dot{w}_k \end{array} \right\} \left\{ \begin{array}{ccc} \dot{u}_k & \dot{v}_k & \dot{w}_k \end{array}\right\} ^{\rm T}~dz \right] dx dy
\label{eqn:kinetic}
\end{equation}
where $\rho_k$ is the mass density of the $k^{th}$ layer, $h_k$ and $h_{k+1}$ are the $z$ coordinates corresponding to the bottom and top surfaces of the $k^{th}$ layer. The strain energy function $U$ is given by
\begin{equation}
U(\delta) = \frac{1}{2} \iint \left[ \sum_{k=1}^n \int\limits_{h_k}^{h_{k+1}} \bvsig^{\rm T} \bveps ~dz \right] dx dy - \iint {\bf q} w~dxdy
\label{eqn:potential}
\end{equation}
where ${\bf q}$ is the distributed force acting on the top surface of the plate. Substituting \Erefs{eqn:potential} and (\ref{eqn:kinetic}) in \Eref{eqn:lagrange}, one obtains the following governing equations for static deflection by neglecting the inertia terms and free vibration of plate.

\paragraph*{Static deflection}
\begin{equation}
\KK \boldsymbol{\delta} = \ff
\label{eqn:staticdefl}
\end{equation}

\paragraph*{Free vibration}
\begin{equation}
\bm \ddot{\boldsymbol{\delta}} + \KK \boldsymbol{\delta} = \bf{0}
\label{eqn:freevib}
\end{equation}
where $\bm$ is the mass matrix, $\KK$ is the stiffness matrix and $\ff$ is the external force vector. In the present study, while performing the integration, terms having thickness co-ordinate $z$ are integrated with higher order Gaussian quadrature, because the material properties vary continuously through the thickness. The terms containing $x$ and $y$ are evaluated using full integration with 3$\times$3 Gauss integration rule. The frequencies and mode shapes are obtained from \Eref{eqn:freevib} using the standard generalized eigenvalue algorithm.

\section{Element description}
\label{eledes}
In this paper, $\mathcal{C}^o$ continuous, eight-noded serendipity quadrilateral shear flexible plate element is used. The finite element represented as per the kinematics based on \Eref{eqn:dispField} is referred to as HSDT13 with cubic variation. The 13 dofs are: $(u_o,v_o,w_o,\theta_x,\theta_y,w_1,\beta_x,\beta_y,\Gamma,\phi_x,\phi_y,\psi_x,\psi_y)$. Five more alternate discrete models are proposed to study the influence of higher-order terms in the displacement functions, whose displacement fields are deduced from the original element by deleting the appropriate degrees of freedom. These alternate models, and the corresponding degrees of freedom are listed in Table \ref{table:alternatemodels}

\begin{table} [htpb]
\renewcommand\arraystretch{1.5}
\caption{Alternate eight-noded finite element models}
\centering
\begin{tabular}{ll}
\hline
Finite element model & Degrees of freedom per node  \\
\hline
HSDT13 & $u_o,v_o,w_o,\theta_x,\theta_y,w_1,\beta_x,\beta_y,\Gamma,\phi_x,\phi_y,\psi_x,\psi_y$ \\
HSDT11A & $u_o,v_o,w_o,\theta_x,\theta_y,\beta_x,\beta_y,\phi_x,\phi_y,\psi_x,\psi_y$ \\
HSDT11B & $u_o,v_o,w_o,\theta_x,\theta_y,w_1,\beta_x,\beta_y,\Gamma,\phi_x,\phi_y$ \\
HSDT9 & $u_o,v_o,w_o,\theta_x,\theta_y,\beta_x,\beta_y,\phi_x,\phi_y$ \\
TSDT7 & $u_o,v_o,w_o,\theta_x,\theta_y,\beta_x,\beta_y$ \\
FSDT5 & $u_o,v_o,w_o,\theta_x,\theta_y$ \\
\hline
\end{tabular}
\label{table:alternatemodels}
\end{table}

\section{Numerical results and discussion}
\label{numexamples}
In this section, we present the static response and the natural frequencies of sandwich plates with homogeneous core and CNT reinforced facesheets using the eight-noded shear flexible quadrilateral element. The effect of plate side-to-thickness ratio, thermal environment, CNT volume fraction on the global response is numerically studied. In this study, only simply supported boundary conditions are considered and are as follows:

\begin{eqnarray}
u_o = w_o = \theta_x = w_1 = \Gamma = \beta_x = \phi_x = \psi_x = 0, \hspace{0.2cm} ~\textup{on} ~ y = 0,b \nonumber \\
v_o = w_o = \theta_y = w_1 = \Gamma = \beta_y = \phi_y = \psi_y = 0, \hspace{0.2cm} ~\textup{on} ~ x= 0,a
\end{eqnarray}
where $a$ and $b$ refer to the length and width of the plate, respectively. The transverse shear stresses are evaluated by integrating the three-dimensional equilibrium equations for all types of elements. For the present study, two different core-to-facesheet thickness $h_H/h_f=$ 2,6 and two thickness ratios $a/h$ (5,10) are considered.

\paragraph{Material properties} For all the numerical studies presented below, unless mentioned otherwise, we employ Poly methyl methacrylate (PMMA) as the matrix in which the CNTs are used as reinforcements. The material properties of which are assumed to be $\rho_m=$ 1150 Kg/m$^3$, $\nu_m=$ 0.34, $\alpha_m=$ 45(1+0.0005$\Delta T$) $\times$ 10$^{-6}$/K and $E_m=$ (3.52-0.0034$T$) GPa in which $T=T_o + \Delta T$ and $T_o=$ 300K. Single walled CNTs are used as reinforcements and the material properties at different temperatures are given in Table \ref{table:cntproperty}. The following CNT efficiency parameters are used~\cite{wangshen2012}: $\eta_1=$ 0.137, $\eta_2=$ 1.002 and $\eta_3=$ 0.715 for $V_{\rm CN}^\ast=$ 0.12; $\eta_1=$ 0.142, $\eta_2=$ 1.626 and $\eta_3=$ 1.138 for $V_{\rm CN}^\ast=$ 0.17 and $\eta_1=$ 0.141, $\eta_2=$ 1.585 and $\eta_3=$ 1.109 for $V_{\rm CN}^\ast=$ 0.28. For this study, we assume that $G_{13} = G_{12}$ and $G_{23}=$ 1.2 $G_{12}$. For the homogeneous core, we use Ti-6Al-4V Titanium alloy. The properties are: $\alpha_H=$ 7.5788(1+6.638$\times$10$^{-4}T$ - 3.147$\times$10$^{-6}T^2$)$\times$10$^{-6}$K, Young's modulus, $E_{_H}=$ 122.56(1-4.586$\times$10$^{-4}T$) GPa, Poisson's ratio $\nu_{_H}=$ 0.29 and mass density $\rho_{_H}=$ 4429 Kg/m$^3$. 

% CNT material properties: Ref: Composites Part B: 43 (2012), 411-421
\begin{table}[htpb]
\centering
\renewcommand\arraystretch{1.5}
\caption{Temperature dependent material properties for $(10,10)$ SWCNT~\cite{wangshen2012}.}
\begin{tabular}{lrrrrr}
\hline
Temperature & $E_{11}^{\rm CN}$ & $E_{22}^{\rm CN}$ & $G_{12}^{\rm CN}$ & $\alpha_{11}^{\rm CN}$ & $\alpha_{22}^{\rm CN}$ \\
$K$ & (TPa) & (TPa) & (TPa) & $(\times$10$^{-6}$/K$)$ & $(\times$10$^{-6}$/K$)$ \\
\hline
300 & 5.6466 & 7.0800 & 1.9445 & 3.4584 & 5.1682 \\
500 & 5.5308 & 6.9348 & 1.9643 & 4.5361 & 5.0189 \\
700 & 5.4744 & 6.8641 & 1.9644 & 4.6677 & 4.8943 \\
\hline
\end{tabular}
\label{table:cntproperty}
\end{table}

\subsection{Static analysis}

% mechanical loading results...
%\input{highOrder_mech}
%\subsubsection{Mechanical loading} 
The static analysis is conducted for carbon nanotube reinforced functionally graded sandwich plate with homogeneous core. The following two types of mechanical loading are considered:

\begin{itemize}
\item Mechanical loading
\begin{itemize}
\item Sinusoidal loading: $q(x,y) = q_o \sin \frac{\pi x}{a} \sin \frac{\pi y}{b}$,
\item Uniformly distributed loading: $q(x,y) = q_o$.
\end{itemize}
\item Thermal loading: $T(x,y) = T_o \left( \frac{2z}{h} \right) \sin \frac{\pi x}{a} \sin \frac{\pi y}{b}$
\end{itemize}
where $q_o$ and $T_o$ are the amplitude of the mechanical load and the thermal load, respectively. The physical quantities are nondimensionalized by relations, unless stated otherwise:
\begin{align}
\renewcommand{\arraystretch}{2}
(\overline{u},\overline{v}) &= \frac{100E_{_H}}{q_o h S^3}(u,v) \nonumber \\
\overline{w} &= \frac{100E_H}{q_o h S^4}w \nonumber \\
\overline{\sigma}_{xx} &= \frac{\sigma_{xx}}{q_oS^2} \nonumber \\
\overline{\sigma}_{xz} &= \frac{\sigma_{xz}}{q_oS}
\end{align}
for the applied mechanical load and by
\begin{align}
\renewcommand{\arraystretch}{2}
(\overline{u},\overline{v}) &= \frac{1}{10 h \alpha_H T_o S}(u,v) \nonumber \\
\overline{w} &= \frac{1}{h \alpha_H T_o h S^4}w\nonumber \\
\overline{\sigma}_{xx} &= \frac{\sigma_{xx}}{100 E_H \alpha_H T_o} \nonumber \\
\overline{\sigma}_{xz} &= \frac{\sigma_{xz}}{10 E_H \alpha_H T_o}
\end{align}
for the applied thermal load, where the $E_{_H}, \alpha_{_H}$ are the Young's modulus and the co-efficient of thermal expansion corresponding to the homogeneous core, evaluated at $T_o=$ 300K and $S=a/h$.

Before proceeding to the detailed analysis of the static response of the sandwich plate for the applied mechanical load, the present formulation is validated against available solutions in the literature.  In this case, only one layer is considered. The CNTs are either uniformly distributed or functionally graded along the thickness direction, given by: 
\begin{equation}
V_{\rm CN}(z)  = \left\{ \begin{array}{cl} V_{\rm CN}^\ast &\textup{UD} \\ \left(1 + \frac{2z}{h} \right) V_{\rm CN}^\ast & \textup{FG-V Type} \\ 2\left(\frac{2|z|}{h} \right) V_{\rm CN}^\ast & \textup{FG-X Type} \end{array} \right.
\end{equation}
The effective material properties, viz., Young's modulus, Poisson's ratio and the mass density are estimated from \Eref{eqn:effecprop}. The effect of type of CNT volume fraction distribution is also considered. PmPV~\cite{zhulei2012} is considered as the matrix with material properties: $E_m=$ 2.1 GPa at room temperature (300K), $\nu_m=$ 0.34. The SWCNTs are chosen as the reinforcements and the material properties for the SWCNT are taken from~\cite{zhulei2012} and given in Table \ref{table:cntproperty}. The CNT efficiency parameters $\eta_j$ are determined according to the effective properties of CNTRCs available by matrching the Young's moduli $E_1$ and $E_2$ with the counterparts computed by rule of mixtures. In this study, the efficiency parameters are: $\eta_1=$ 0.149 and $\eta_2=$ 0.934 for the case of $V_{\rm CN}^\ast=$ 0.11; $\eta_1=$ 0.150 and $\eta_2=$ 0.941 for the case of $V_{\rm CN}^\ast=$ 0.14 and $\eta_1=$ 0.149 and $\eta_2=$ 1.381 for the case of $V_{\rm CN}^\ast=$ 0.17. It is assumed that $\eta_2=\eta_3$ for this study.

% mechanical validation with results from Ref: Composite Structures, v94, 2012, 1450-1460
\begin{table}[htpb]
\centering
\renewcommand\arraystretch{1.5}
\caption{Convergence of maximum center deflection $w_c = -w_o/h \times$10$^{-2}$ with mesh size for a simply supported square plate subjected to uniformly distributed load with $a/h=$ 20. The volume fraction of the CNT is $V_{\rm CN}^\ast=$ 0.14.}
\begin{tabular}{lrrrrrrr}
\hline
Mesh & \multicolumn{3}{c}{FSDT5} && \multicolumn{3}{c}{HSDT11B}\\
\cline{2-4}\cline{6-8}
& UD & FG-V & FG-X &&  UD & FG-V & FG-X\\
\hline
4$\times$4 & 3.0174 & 4.0418 & 2.2710 && 3.0172 & 4.0959 & 2.2764\\
6$\times$6 & 3.0018 & 4.0209 & 2.2593 && 3.0017 &  4.0748 & 2.2647 \\
8$\times$8 & 2.9993 & 4.0176 & 2.2574 && 2.9873 & 4.0552 & 2.2538\\
16$\times$16 & 2.9993 & 4.0176 & 2.2574 && 2.9873 & 4.0552 & 2.2538\\
Ref.~\cite{zhulei2012} & 3.0010 & 4.0250 & 2.2560 && - & - & - \\
\hline
\end{tabular}
\label{table:mechvalid}
\end{table}

% mechanical validation with results from Ref: Composite Structures, v94, 2012, 1450-1460
\begin{table}[htpb]
\centering
\renewcommand\arraystretch{1.5}
\caption{Convergence of maximum center deflection $w_c = -w_o/h \times$10$^{-2}$ with mesh size for a simply supported square sandwich plate subjected to sinusoidally distributed mechanical load with $a/h=$ 5.}
\begin{tabular}{lrr}
\hline
Mesh & \multicolumn{2}{c}{$h_H/h_f=$ 2}\\
\cline{2-3}
& $V_{\rm CN}^\ast=$ 0.17 & $V_{\rm CN}^\ast=$ 0.28\\
\hline
4$\times$4 & 0.1042 & 0.0911\\
6$\times$6 &   0.1037 & 0.0907\\
8$\times$8 &  0.1036 & 0.0906\\
16$\times$16 &  0.1036 & 0.0906\\
\hline
\end{tabular}
\label{table:mechvalidsandwich}
\end{table}

Tables \ref{table:mechvalid} - \ref{table:mechvalidsandwich} present the convergence of the maximum center deflection with decreasing mesh size for a simply supported square plate subjected to uniformly distributed load with $a/h=$ 20 employing the first- (FSDT5) \& higher order (HSDT11B) and with $a/h=$ 5 employing the higher-order (HSDT13) structural models. The results of the first-order shear deformation theory (FSDT5) employed here is compared with the available results from the literature~\cite{zhulei2012}. It can be seen that the results from the present formulation are in excellent agreement with the results in the literature and that a 8 $\times$ 8 mesh is found to be adequate to model the full plate, viz., sandwich plate or single layer for the present analysis. The influence of the CNT volume fraction distribution is also studied. It is seen that the plate with CNT reinforced in FG-V has the maximum center deflection and the plate with FG-X type of reinforcement has the minimum center deflection.

Next the numerical study is carried out to study the influence of the plate thickness $a/h$, the core-to-facesheet thickness $h_H/h_f$ and the CNT volume fraction $V_{\rm CN}^\ast$ on the deflection and the stresses for a plate with sinusoidal loading. Numerical results are tabulated in Tables \ref{table:vcnt17influence} and \ref{table:vcnt28influence} for an applied sinusoidal mechanical load and in Tables \ref{table:vcnt17influencetherm} and \ref{table:vcnt28influencetherm} for an applied thermal load. It is inferred from these Tables that with increasing the CNT volume fraction $V_{\rm CN}^\ast$, the non-dimensionalized displacement decreases whereas the change in the non-dimensionalized stresses however depends on the thickness ratio. Furthermore, the non-dimensionalized displacements and the non-dimensionalized stresses decrease in general with increasing core thickness. This can be attributed to the change in the flexural stiffness due to increase in the core thickness and the addition of CNTs to the facesheets. It can be also noted that the non-dimensionalized displacements and the non-dimensionalized stresses predicted at the neutral surface using HSDT13, HSDT11A and HSDT11B are different from those of HSDT9, TSDT7 and FSDT5. However, the latter models cannot predict the through thickness stretching effects in the displacements and the stresses. 

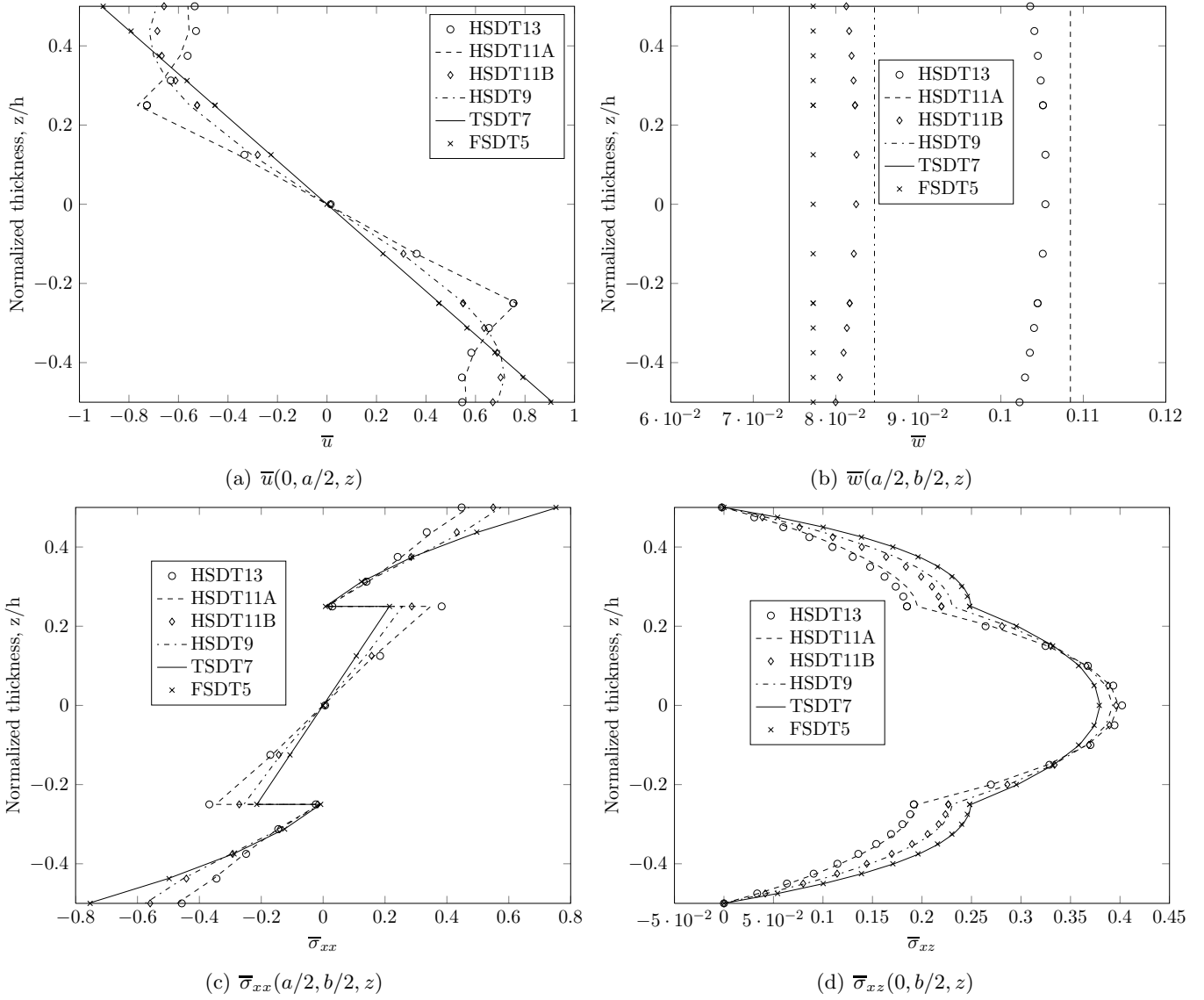
\begin{figure}[htpb]
\newlength\figureheight 
\newlength\figurewidth 
\setlength\figureheight{8cm} 
\setlength\figurewidth{10cm}
\centering
\subfigure[$\overline{u}(0,a/2,z)$]{\scalebox{0.75}{% This file was created by matlab2tikz v0.3.3.
% Copyright (c) 2008--2013, Nico Schlömer <nico.schloemer@gmail.com>
% All rights reserved.
% 
% The latest updates can be retrieved from
%   http://www.mathworks.com/matlabcentral/fileexchange/22022-matlab2tikz
% where you can also make suggestions and rate matlab2tikz.
% 
% 
% 
\begin{tikzpicture}

\begin{axis}[%
width=\figurewidth,
height=\figureheight,
scale only axis,
xmin=-1,
xmax=1,
xlabel={$\overline{u}$},
ymin=-0.5,
ymax=0.5,
ylabel={Normalized thickness, z/h},
legend style={draw=black,fill=white,legend cell align=left}
]
\addplot [
color=black,
only marks,
mark=o,
mark options={solid}
]
table[row sep=crcr]{
0.546282 -0.5\\
0.54517 -0.4375\\
0.582964 -0.375\\
0.654066 -0.3125\\
0.752877 -0.25\\
0.752877 -0.25\\
0.362104 -0.125\\
0.0149878 0\\
-0.33326 0.125\\
-0.727426 0.25\\
-0.727426 0.25\\
-0.63116 0.3125\\
-0.56317 0.375\\
-0.529052 0.4375\\
-0.534406 0.5\\
};
\addlegendentry{HSDT13};

\addplot [
color=black,
dashed
]
table[row sep=crcr]{
0.561744 -0.5\\
0.557207 -0.4375\\
0.593551 -0.375\\
0.664936 -0.3125\\
0.765522 -0.25\\
0.765522 -0.25\\
0.3594 -0.125\\
0 0\\
-0.3594 0.125\\
-0.765522 0.25\\
-0.765522 0.25\\
-0.664936 0.3125\\
-0.593551 0.375\\
-0.557207 0.4375\\
-0.561744 0.5\\
};
\addlegendentry{HSDT11A};

\addplot [
color=black,
only marks,
mark=diamond,
mark options={solid}
]
table[row sep=crcr]{
0.670051 -0.5\\
0.701335 -0.4375\\
0.687465 -0.375\\
0.634852 -0.3125\\
0.549906 -0.25\\
0.549906 -0.25\\
0.308653 -0.125\\
0.0149878 0\\
-0.279809 0.125\\
-0.524455 0.25\\
-0.524455 0.25\\
-0.611947 0.3125\\
-0.667671 0.375\\
-0.685217 0.4375\\
-0.658175 0.5\\
};
\addlegendentry{HSDT11B};

\addplot [
color=black,
dash pattern=on 1pt off 3pt on 3pt off 3pt
]
table[row sep=crcr]{
0.688301 -0.5\\
0.717563 -0.4375\\
0.700705 -0.375\\
0.644315 -0.3125\\
0.554984 -0.25\\
0.554984 -0.25\\
0.303846 -0.125\\
0 0\\
-0.303846 0.125\\
-0.554984 0.25\\
-0.554984 0.25\\
-0.644315 0.3125\\
-0.700705 0.375\\
-0.717563 0.4375\\
-0.688301 0.5\\
};
\addlegendentry{HSDT9};

\addplot [
color=black,
solid
]
table[row sep=crcr]{
0.910794 -0.5\\
0.796945 -0.4375\\
0.683096 -0.375\\
0.569246 -0.3125\\
0.455397 -0.25\\
0.455397 -0.25\\
0.227699 -0.125\\
0 0\\
-0.227699 0.125\\
-0.455397 0.25\\
-0.455397 0.25\\
-0.569246 0.3125\\
-0.683096 0.375\\
-0.796945 0.4375\\
-0.910794 0.5\\
};
\addlegendentry{TSDT7};

\addplot [
color=black,
only marks,
mark=x,
mark options={solid}
]
table[row sep=crcr]{
0.904212 -0.5\\
0.791186 -0.4375\\
0.678159 -0.375\\
0.565133 -0.3125\\
0.452106 -0.25\\
0.452106 -0.25\\
0.226053 -0.125\\
0 0\\
-0.226053 0.125\\
-0.452106 0.25\\
-0.452106 0.25\\
-0.565133 0.3125\\
-0.678159 0.375\\
-0.791186 0.4375\\
-0.904212 0.5\\
};
\addlegendentry{FSDT5};

\end{axis}
\end{tikzpicture}%
}}
\subfigure[$\overline{w}(a/2,b/2,z)$]{\scalebox{0.75}{% This file was created by matlab2tikz v0.3.3.
% Copyright (c) 2008--2013, Nico Schlömer <nico.schloemer@gmail.com>
% All rights reserved.
% 
% The latest updates can be retrieved from
%   http://www.mathworks.com/matlabcentral/fileexchange/22022-matlab2tikz
% where you can also make suggestions and rate matlab2tikz.
% 
% 
% 
\begin{tikzpicture}

\begin{axis}[%
width=\figurewidth,
height=\figureheight,
scale only axis,
xmin=0.06,
xmax=0.12,
xlabel={$\overline{w}$},
ymin=-0.5,
ymax=0.5,
ylabel={Normalized thickness, z/h},
legend style={at={(0.42,0.5)},anchor=south west,draw=black,fill=white,legend cell align=left}
]
\addplot [
color=black,
only marks,
mark=o,
mark options={solid}
]
table[row sep=crcr]{
0.102274 -0.5\\
0.102938 -0.4375\\
0.103524 -0.375\\
0.104032 -0.3125\\
0.104463 -0.25\\
0.104463 -0.25\\
0.105091 -0.125\\
0.105408 0\\
0.105414 0.125\\
0.10511 0.25\\
0.10511 0.25\\
0.104841 0.3125\\
0.104494 0.375\\
0.10407 0.4375\\
0.103567 0.5\\
};
\addlegendentry{HSDT13};

\addplot [
color=black,
dashed
]
table[row sep=crcr]{
0.108439 -0.5\\
0.108439 -0.4375\\
0.108439 -0.375\\
0.108439 -0.3125\\
0.108439 -0.25\\
0.108439 -0.25\\
0.108439 -0.125\\
0.108439 0\\
0.108439 0.125\\
0.108439 0.25\\
0.108439 0.25\\
0.108439 0.3125\\
0.108439 0.375\\
0.108439 0.4375\\
0.108439 0.5\\
};
\addlegendentry{HSDT11A};

\addplot [
color=black,
only marks,
mark=diamond,
mark options={solid}
]
table[row sep=crcr]{
0.0799828 -0.5\\
0.0804948 -0.4375\\
0.0809493 -0.375\\
0.0813463 -0.3125\\
0.0816859 -0.25\\
0.0816859 -0.25\\
0.0821925 -0.125\\
0.0824691 0\\
0.0825159 0.125\\
0.0823326 0.25\\
0.0823326 0.25\\
0.0821548 0.3125\\
0.0819194 0.375\\
0.0816266 0.4375\\
0.0812763 0.5\\
};
\addlegendentry{HSDT11B};

\addplot [
color=black,
dash pattern=on 1pt off 3pt on 3pt off 3pt
]
table[row sep=crcr]{
0.0846917 -0.5\\
0.0846917 -0.4375\\
0.0846917 -0.375\\
0.0846917 -0.3125\\
0.0846917 -0.25\\
0.0846917 -0.25\\
0.0846917 -0.125\\
0.0846917 0\\
0.0846917 0.125\\
0.0846917 0.25\\
0.0846917 0.25\\
0.0846917 0.3125\\
0.0846917 0.375\\
0.0846917 0.4375\\
0.0846917 0.5\\
};
\addlegendentry{HSDT9};

\addplot [
color=black,
solid
]
table[row sep=crcr]{
0.0743557 -0.5\\
0.0743557 -0.4375\\
0.0743557 -0.375\\
0.0743557 -0.3125\\
0.0743557 -0.25\\
0.0743557 -0.25\\
0.0743557 -0.125\\
0.0743557 0\\
0.0743557 0.125\\
0.0743557 0.25\\
0.0743557 0.25\\
0.0743557 0.3125\\
0.0743557 0.375\\
0.0743557 0.4375\\
0.0743557 0.5\\
};
\addlegendentry{TSDT7};

\addplot [
color=black,
only marks,
mark=x,
mark options={solid}
]
table[row sep=crcr]{
0.0772578 -0.5\\
0.0772578 -0.4375\\
0.0772578 -0.375\\
0.0772578 -0.3125\\
0.0772578 -0.25\\
0.0772578 -0.25\\
0.0772578 -0.125\\
0.0772578 0\\
0.0772578 0.125\\
0.0772578 0.25\\
0.0772578 0.25\\
0.0772578 0.3125\\
0.0772578 0.375\\
0.0772578 0.4375\\
0.0772578 0.5\\
};
\addlegendentry{FSDT5};

\end{axis}
\end{tikzpicture}%
}}
\subfigure[$\overline{\sigma}_{xx}(a/2,b/2,z)$]{\scalebox{0.75}{% This file was created by matlab2tikz v0.3.3.
% Copyright (c) 2008--2013, Nico Schlömer <nico.schloemer@gmail.com>
% All rights reserved.
% 
% The latest updates can be retrieved from
%   http://www.mathworks.com/matlabcentral/fileexchange/22022-matlab2tikz
% where you can also make suggestions and rate matlab2tikz.
% 
% 
% 
\begin{tikzpicture}

\begin{axis}[%
width=\figurewidth,
height=\figureheight,
scale only axis,
xmin=-0.8,
xmax=0.8,
xlabel={$\overline{\sigma}_{xx}$},
ymin=-0.5,
ymax=0.5,
ylabel={Normalized thickness, z/h},
legend style={at={(0.152678571428571,0.5)},anchor=south west,draw=black,fill=white,legend cell align=left}
]
\addplot [
color=black,
only marks,
mark=o,
mark options={solid}
]
table[row sep=crcr]{
-0.456845 -0.5\\
-0.344679 -0.4375\\
-0.248945 -0.375\\
-0.145095 -0.3125\\
-0.0238522 -0.25\\
-0.367914 -0.25\\
-0.170669 -0.125\\
0.00642642 0\\
0.184096 0.125\\
0.383063 0.25\\
0.0299892 0.25\\
0.140319 0.3125\\
0.240863 0.375\\
0.33486 0.4375\\
0.447572 0.5\\
};
\addlegendentry{HSDT13};

\addplot [
color=black,
dashed
]
table[row sep=crcr]{
-0.471868 -0.5\\
-0.353621 -0.4375\\
-0.253981 -0.375\\
-0.147055 -0.3125\\
-0.0125063 -0.25\\
-0.347551 -0.25\\
-0.164987 -0.125\\
0 0\\
0.164987 0.125\\
0.347551 0.25\\
0.0121206 0.25\\
0.146121 0.3125\\
0.253147 0.375\\
0.352838 0.4375\\
0.471473 0.5\\
};
\addlegendentry{HSDT11A};

\addplot [
color=black,
only marks,
mark=diamond,
mark options={solid}
]
table[row sep=crcr]{
-0.558994 -0.5\\
-0.441944 -0.4375\\
-0.292586 -0.375\\
-0.140477 -0.3125\\
-0.0177082 -0.25\\
-0.271157 -0.25\\
-0.143234 -0.125\\
0.00642642 0\\
0.156661 0.125\\
0.286306 0.25\\
0.0230672 0.25\\
0.135727 0.3125\\
0.284357 0.375\\
0.431907 0.4375\\
0.549634 0.5\\
};
\addlegendentry{HSDT11B};

\addplot [
color=black,
dash pattern=on 1pt off 3pt on 3pt off 3pt
]
table[row sep=crcr]{
-0.57518 -0.5\\
-0.45236 -0.4375\\
-0.297817 -0.375\\
-0.141731 -0.3125\\
-0.00931051 -0.25\\
-0.256301 -0.25\\
-0.137207 -0.125\\
0 0\\
0.137207 0.125\\
0.256301 0.25\\
0.00903678 0.25\\
0.140826 0.3125\\
0.296832 0.375\\
0.451351 0.4375\\
0.574696 0.5\\
};
\addlegendentry{HSDT9};

\addplot [
color=black,
solid
]
table[row sep=crcr]{
-0.758434 -0.5\\
-0.501202 -0.4375\\
-0.289912 -0.375\\
-0.125203 -0.3125\\
-0.00785136 -0.25\\
-0.214066 -0.25\\
-0.107033 -0.125\\
0 0\\
0.107033 0.125\\
0.214066 0.25\\
0.00763189 0.25\\
0.124403 0.3125\\
0.288953 0.375\\
0.500082 0.4375\\
0.757794 0.5\\
};
\addlegendentry{TSDT7};

\addplot [
color=black,
only marks,
mark=x,
mark options={solid}
]
table[row sep=crcr]{
-0.753141 -0.5\\
-0.497742 -0.4375\\
-0.287955 -0.375\\
-0.124412 -0.3125\\
-0.0079077 -0.25\\
-0.214525 -0.25\\
-0.107263 -0.125\\
0 0\\
0.107263 0.125\\
0.214525 0.25\\
0.00769259 0.25\\
0.123619 0.3125\\
0.287002 0.375\\
0.496631 0.4375\\
0.752505 0.5\\
};
\addlegendentry{FSDT5};

\end{axis}
\end{tikzpicture}%
}}
\subfigure[$\overline{\sigma}_{xz}(0,b/2,z)$]{\scalebox{0.75}{% This file was created by matlab2tikz v0.3.3.
% Copyright (c) 2008--2013, Nico Schlömer <nico.schloemer@gmail.com>
% All rights reserved.
% 
% The latest updates can be retrieved from
%   http://www.mathworks.com/matlabcentral/fileexchange/22022-matlab2tikz
% where you can also make suggestions and rate matlab2tikz.
% 
% 
% 
\begin{tikzpicture}

\begin{axis}[%
width=\figurewidth,
height=\figureheight,
scale only axis,
xmin=-0.05,
xmax=0.45,
xlabel={$\overline{\sigma}_{xz}$},
ymin=-0.5,
ymax=0.5,
ylabel={Normalized thickness, z/h},
legend style={at={(0.152678571428571,0.4)},anchor=south west,draw=black,fill=white,legend cell align=left}
]
\addplot [
color=black,
only marks,
mark=o,
mark options={solid}
]
table[row sep=crcr]{
0 -0.5\\
0.0335423 -0.475\\
0.0636466 -0.45\\
0.0906128 -0.425\\
0.114612 -0.4\\
0.135693 -0.375\\
0.153794 -0.35\\
0.168748 -0.325\\
0.180294 -0.3\\
0.188088 -0.275\\
0.191788 -0.25\\
0.191788 -0.25\\
0.269571 -0.2\\
0.328696 -0.15\\
0.370094 -0.1\\
0.394424 -0.05\\
0.402076 0\\
0.393169 0.05\\
0.367549 0.1\\
0.324793 0.15\\
0.264207 0.2\\
0.184826 0.25\\
0.184826 0.25\\
0.18112 0.275\\
0.173474 0.3\\
0.162204 0.325\\
0.147642 0.35\\
0.130038 0.375\\
0.109545 0.4\\
0.0862105 0.425\\
0.0599696 0.45\\
0.030635 0.475\\
-0.00211132 0.5\\
};
\addlegendentry{HSDT13};

\addplot [
color=black,
dashed
]
table[row sep=crcr]{
0 -0.5\\
0.0346483 -0.475\\
0.065601 -0.45\\
0.0932112 -0.425\\
0.117693 -0.4\\
0.139128 -0.375\\
0.15748 -0.35\\
0.172597 -0.325\\
0.184226 -0.3\\
0.192021 -0.275\\
0.195562 -0.25\\
0.195562 -0.25\\
0.270108 -0.2\\
0.326931 -0.15\\
0.3669 -0.1\\
0.390633 -0.05\\
0.398503 0\\
0.390633 0.05\\
0.3669 0.1\\
0.326931 0.15\\
0.270108 0.2\\
0.195562 0.25\\
0.195562 0.25\\
0.192021 0.275\\
0.184226 0.3\\
0.172597 0.325\\
0.15748 0.35\\
0.139128 0.375\\
0.117693 0.4\\
0.0932112 0.425\\
0.065601 0.45\\
0.0346483 0.475\\
0 0.5\\
};
\addlegendentry{HSDT11A};

\addplot [
color=black,
only marks,
mark=diamond,
mark options={solid}
]
table[row sep=crcr]{
0 -0.5\\
0.0417192 -0.475\\
0.0800143 -0.45\\
0.114315 -0.425\\
0.144214 -0.4\\
0.169461 -0.375\\
0.189947 -0.35\\
0.205695 -0.325\\
0.216854 -0.3\\
0.223684 -0.275\\
0.226611 -0.25\\
0.226611 -0.25\\
0.286228 -0.2\\
0.333898 -0.15\\
0.368545 -0.1\\
0.389398 -0.05\\
0.395987 0\\
0.388143 0.05\\
0.366 0.1\\
0.329995 0.15\\
0.280865 0.2\\
0.21965 0.25\\
0.21965 0.25\\
0.216717 0.275\\
0.210034 0.3\\
0.199151 0.325\\
0.183796 0.35\\
0.163807 0.375\\
0.139148 0.4\\
0.109912 0.425\\
0.0763373 0.45\\
0.0388118 0.475\\
-0.00211132 0.5\\
};
\addlegendentry{HSDT11B};

\addplot [
color=black,
dash pattern=on 1pt off 3pt on 3pt off 3pt
]
table[row sep=crcr]{
0 -0.5\\
0.0429387 -0.475\\
0.0822023 -0.45\\
0.117248 -0.425\\
0.147698 -0.4\\
0.173326 -0.375\\
0.194049 -0.35\\
0.209919 -0.325\\
0.221106 -0.3\\
0.227893 -0.275\\
0.230675 -0.25\\
0.230675 -0.25\\
0.287363 -0.2\\
0.332709 -0.15\\
0.365772 -0.1\\
0.385879 -0.05\\
0.392626 0\\
0.385879 0.05\\
0.365772 0.1\\
0.332709 0.15\\
0.287363 0.2\\
0.230675 0.25\\
0.230675 0.25\\
0.227893 0.275\\
0.221106 0.3\\
0.209919 0.325\\
0.194049 0.35\\
0.173326 0.375\\
0.147698 0.4\\
0.117248 0.425\\
0.0822023 0.45\\
0.0429387 0.475\\
0 0.5\\
};
\addlegendentry{HSDT9};

\addplot [
color=black,
solid
]
table[row sep=crcr]{
0 -0.5\\
0.0545021 -0.475\\
0.100909 -0.45\\
0.139798 -0.425\\
0.171742 -0.4\\
0.197318 -0.375\\
0.217099 -0.35\\
0.231657 -0.325\\
0.241568 -0.3\\
0.247403 -0.275\\
0.249744 -0.25\\
0.249744 -0.25\\
0.296575 -0.2\\
0.332999 -0.15\\
0.359016 -0.1\\
0.374626 -0.05\\
0.37983 0\\
0.374626 0.05\\
0.359016 0.1\\
0.332999 0.15\\
0.296575 0.2\\
0.249744 0.25\\
0.249744 0.25\\
0.247403 0.275\\
0.241568 0.3\\
0.231657 0.325\\
0.217099 0.35\\
0.197318 0.375\\
0.171742 0.4\\
0.139798 0.425\\
0.100909 0.45\\
0.0545021 0.475\\
0 0.5\\
};
\addlegendentry{TSDT7};

\addplot [
color=black,
only marks,
mark=x,
mark options={solid}
]
table[row sep=crcr]{
0 -0.5\\
0.054137 -0.475\\
0.100235 -0.45\\
0.138867 -0.425\\
0.170604 -0.4\\
0.196015 -0.375\\
0.215672 -0.35\\
0.230143 -0.325\\
0.239998 -0.3\\
0.245805 -0.275\\
0.248143 -0.25\\
0.248143 -0.25\\
0.295272 -0.2\\
0.331928 -0.15\\
0.35811 -0.1\\
0.37382 -0.05\\
0.379057 0\\
0.37382 0.05\\
0.35811 0.1\\
0.331928 0.15\\
0.295272 0.2\\
0.248143 0.25\\
0.248143 0.25\\
0.245805 0.275\\
0.239998 0.3\\
0.230143 0.325\\
0.215672 0.35\\
0.196015 0.375\\
0.170604 0.4\\
0.138867 0.425\\
0.100235 0.45\\
0.054137 0.475\\
0 0.5\\
};
\addlegendentry{FSDT5};

\end{axis}
\end{tikzpicture}%
}}
\caption{Displacements and stresses through the thickness for the square plates with simply supported edges with $a/h=$ 5, $h_H/h_f=$ 2 and $V_{\rm CN}^\ast=$ 0.17 and the CNT distribution is FG-X. The plate is subjected to a mechanical load $q(x,y) = q_o\sin (\pi x/a) \sin (\pi y/b)$}
\label{fig:thickdispstressmech}
\end{figure}

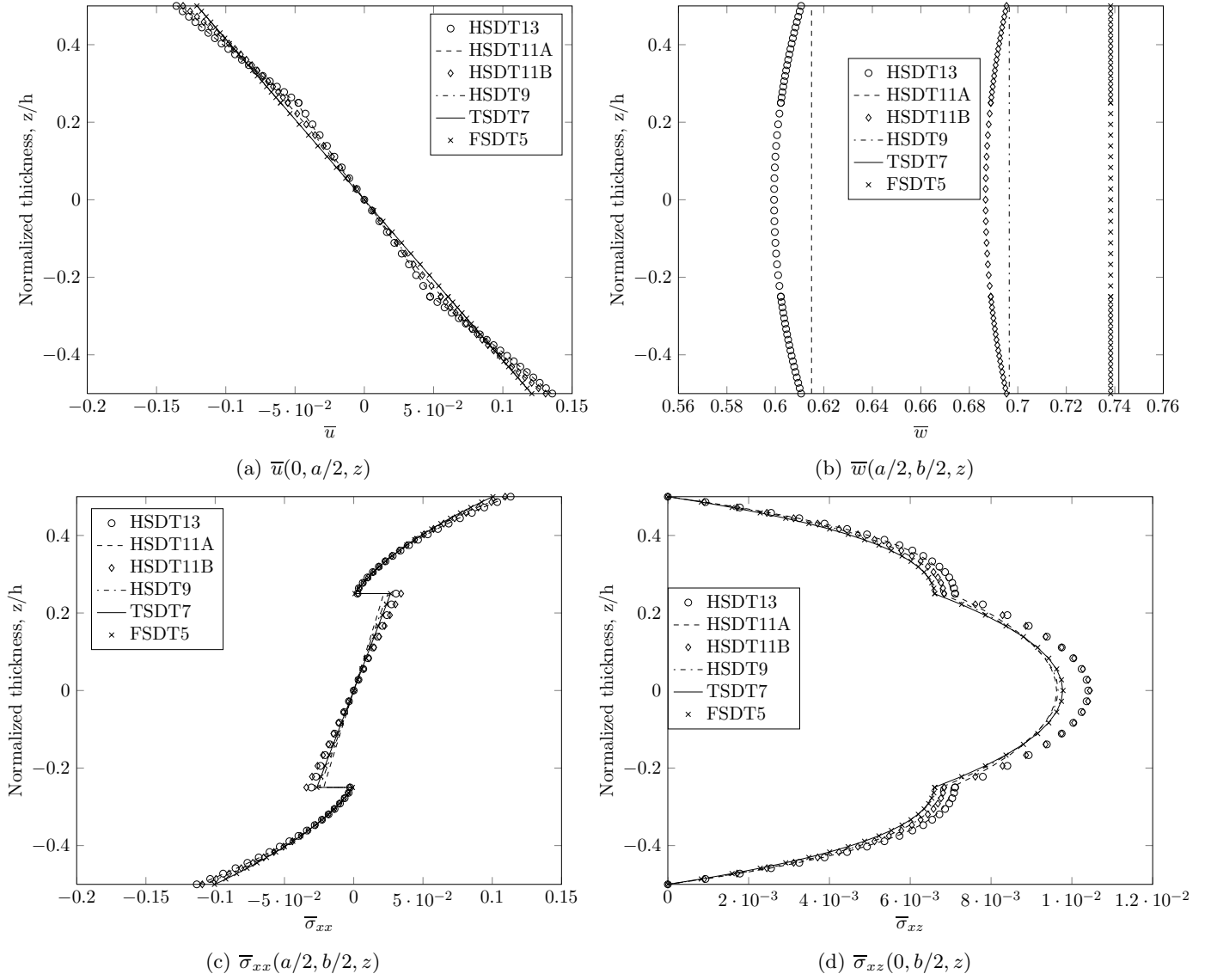
\begin{figure}[htpb]
%\newlength\figureheight 
%\newlength\figurewidth 
\setlength\figureheight{8cm} 
\setlength\figurewidth{10cm}
\centering
\subfigure[$\overline{u}(0,a/2,z)$]{\scalebox{0.75}{% This file was created by matlab2tikz v0.3.3.
% Copyright (c) 2008--2013, Nico Schlömer <nico.schloemer@gmail.com>
% All rights reserved.
% 
% The latest updates can be retrieved from
%   http://www.mathworks.com/matlabcentral/fileexchange/22022-matlab2tikz
% where you can also make suggestions and rate matlab2tikz.
% 
% 
% 
\begin{tikzpicture}

\begin{axis}[%
width=\figurewidth,
height=\figureheight,
scale only axis,
xmin=-0.2,
xmax=0.15,
xlabel={$\overline{u}$},
ymin=-0.5,
ymax=0.5,
ylabel={Normalized thickness, z/h},
legend style={draw=black,fill=white,legend cell align=left}
]
\addplot [
color=black,
only marks,
mark=o,
mark options={solid}
]
table[row sep=crcr]{
0.135674 -0.5\\
0.131167 -0.486111\\
0.126606 -0.472222\\
0.121993 -0.4583335\\
0.117328 -0.4444445\\
0.112614 -0.4305555\\
0.107851 -0.4166665\\
0.103042 -0.402778\\
0.0981873 -0.388889\\
0.0932895 -0.375\\
0.0883498 -0.361111\\
0.0833697 -0.347222\\
0.0783508 -0.3333335\\
0.0732946 -0.3194445\\
0.0682027 -0.3055555\\
0.0630766 -0.2916665\\
0.057918 -0.277778\\
0.0527282 -0.263889\\
0.0475089 -0.25\\
0.0475089 -0.25\\
0.0425118 -0.222222\\
0.0374152 -0.1944445\\
0.0322317 -0.1666665\\
0.0269737 -0.138889\\
0.0216535 -0.111111\\
0.0162836 -0.0833335\\
0.0108764 -0.0555555\\
0.00544443 -0.0277778\\
0 0\\
-0.00544443 0.0277778\\
-0.0108764 0.0555555\\
-0.0162836 0.0833335\\
-0.0216535 0.111111\\
-0.0269737 0.138889\\
-0.0322317 0.1666665\\
-0.0374152 0.1944445\\
-0.0425118 0.222222\\
-0.0475089 0.25\\
-0.0475089 0.25\\
-0.0527282 0.263889\\
-0.057918 0.277778\\
-0.0630766 0.2916665\\
-0.0682027 0.3055555\\
-0.0732946 0.3194445\\
-0.0783508 0.3333335\\
-0.0833697 0.347222\\
-0.0883498 0.361111\\
-0.0932895 0.375\\
-0.0981873 0.388889\\
-0.103042 0.402778\\
-0.107851 0.4166665\\
-0.112614 0.4305555\\
-0.117328 0.4444445\\
-0.121993 0.4583335\\
-0.126606 0.472222\\
-0.131167 0.486111\\
-0.135674 0.5\\
};
\addlegendentry{HSDT13};

\addplot [
color=black,
dashed
]
table[row sep=crcr]{
0.134973 -0.5\\
0.130573 -0.486111\\
0.126119 -0.472222\\
0.121612 -0.4583335\\
0.117054 -0.4444445\\
0.112445 -0.4305555\\
0.107789 -0.4166665\\
0.103085 -0.402778\\
0.0983368 -0.388889\\
0.0935447 -0.375\\
0.0887105 -0.361111\\
0.0838359 -0.347222\\
0.0789224 -0.3333335\\
0.0739715 -0.3194445\\
0.0689848 -0.3055555\\
0.0639638 -0.2916665\\
0.0589101 -0.277778\\
0.0538253 -0.263889\\
0.0487109 -0.25\\
0.0487109 -0.25\\
0.0435809 -0.222222\\
0.0383513 -0.1944445\\
0.0330345 -0.1666665\\
0.0276429 -0.138889\\
0.0221891 -0.111111\\
0.0166854 -0.0833335\\
0.0111444 -0.0555555\\
0.00557841 -0.0277778\\
0 0\\
-0.00557841 0.0277778\\
-0.0111444 0.0555555\\
-0.0166854 0.0833335\\
-0.0221891 0.111111\\
-0.0276429 0.138889\\
-0.0330345 0.1666665\\
-0.0383513 0.1944445\\
-0.0435809 0.222222\\
-0.0487109 0.25\\
-0.0487109 0.25\\
-0.0538253 0.263889\\
-0.0589101 0.277778\\
-0.0639638 0.2916665\\
-0.0689848 0.3055555\\
-0.0739715 0.3194445\\
-0.0789224 0.3333335\\
-0.0838359 0.347222\\
-0.0887105 0.361111\\
-0.0935447 0.375\\
-0.0983368 0.388889\\
-0.103085 0.402778\\
-0.107789 0.4166665\\
-0.112445 0.4305555\\
-0.117054 0.4444445\\
-0.121612 0.4583335\\
-0.126119 0.472222\\
-0.130573 0.486111\\
-0.134973 0.5\\
};
\addlegendentry{HSDT11A};

\addplot [
color=black,
only marks,
mark=diamond,
mark options={solid}
]
table[row sep=crcr]{
0.130842 -0.5\\
0.125761 -0.486111\\
0.120802 -0.472222\\
0.115962 -0.4583335\\
0.111237 -0.4444445\\
0.106624 -0.4305555\\
0.102119 -0.4166665\\
0.097719 -0.402778\\
0.0934203 -0.388889\\
0.0892194 -0.375\\
0.0851128 -0.361111\\
0.0810971 -0.347222\\
0.0771686 -0.3333335\\
0.0733239 -0.3194445\\
0.0695596 -0.3055555\\
0.0658722 -0.2916665\\
0.0622581 -0.277778\\
0.0587138 -0.263889\\
0.0552359 -0.25\\
0.0552359 -0.25\\
0.0484652 -0.222222\\
0.0419181 -0.1944445\\
0.0355665 -0.1666665\\
0.0293826 -0.138889\\
0.0233385 -0.111111\\
0.017406 -0.0833335\\
0.0115575 -0.0555555\\
0.00576476 -0.0277778\\
0 0\\
-0.00576476 0.0277778\\
-0.0115575 0.0555555\\
-0.017406 0.0833335\\
-0.0233385 0.111111\\
-0.0293826 0.138889\\
-0.0355665 0.1666665\\
-0.0419181 0.1944445\\
-0.0484652 0.222222\\
-0.0552359 0.25\\
-0.0552359 0.25\\
-0.0587138 0.263889\\
-0.0622581 0.277778\\
-0.0658722 0.2916665\\
-0.0695596 0.3055555\\
-0.0733239 0.3194445\\
-0.0771686 0.3333335\\
-0.0810971 0.347222\\
-0.0851128 0.361111\\
-0.0892194 0.375\\
-0.0934203 0.388889\\
-0.097719 0.402778\\
-0.102119 0.4166665\\
-0.106624 0.4305555\\
-0.111237 0.4444445\\
-0.115962 0.4583335\\
-0.120802 0.472222\\
-0.125761 0.486111\\
-0.130842 0.5\\
};
\addlegendentry{HSDT11B};

\addplot [
color=black,
dash pattern=on 1pt off 3pt on 3pt off 3pt
]
table[row sep=crcr]{
0.13046 -0.5\\
0.125518 -0.486111\\
0.120688 -0.472222\\
0.115966 -0.4583335\\
0.11135 -0.4444445\\
0.106834 -0.4305555\\
0.102418 -0.4166665\\
0.0980969 -0.402778\\
0.0938681 -0.388889\\
0.0897283 -0.375\\
0.0856745 -0.361111\\
0.0817033 -0.347222\\
0.0778117 -0.3333335\\
0.0739964 -0.3194445\\
0.0702543 -0.3055555\\
0.0665822 -0.2916665\\
0.0629768 -0.277778\\
0.0594351 -0.263889\\
0.0559539 -0.25\\
0.0559539 -0.25\\
0.0491599 -0.222222\\
0.0425696 -0.1944445\\
0.0361574 -0.1666665\\
0.0298979 -0.138889\\
0.0237656 -0.111111\\
0.0177351 -0.0833335\\
0.011781 -0.0555555\\
0.00587778 -0.0277778\\
0 0\\
-0.00587778 0.0277778\\
-0.011781 0.0555555\\
-0.0177351 0.0833335\\
-0.0237656 0.111111\\
-0.0298979 0.138889\\
-0.0361574 0.1666665\\
-0.0425696 0.1944445\\
-0.0491599 0.222222\\
-0.0559539 0.25\\
-0.0559539 0.25\\
-0.0594351 0.263889\\
-0.0629768 0.277778\\
-0.0665822 0.2916665\\
-0.0702543 0.3055555\\
-0.0739964 0.3194445\\
-0.0778117 0.3333335\\
-0.0817033 0.347222\\
-0.0856745 0.361111\\
-0.0897283 0.375\\
-0.0938681 0.388889\\
-0.0980969 0.402778\\
-0.102418 0.4166665\\
-0.106834 0.4305555\\
-0.11135 0.4444445\\
-0.115966 0.4583335\\
-0.120688 0.472222\\
-0.125518 0.486111\\
-0.13046 0.5\\
};
\addlegendentry{HSDT9};

\addplot [
color=black,
solid
]
table[row sep=crcr]{
0.12063 -0.5\\
0.117279 -0.486111\\
0.113928 -0.472222\\
0.110577 -0.4583335\\
0.107227 -0.4444445\\
0.103876 -0.4305555\\
0.100525 -0.4166665\\
0.097174 -0.402778\\
0.0938232 -0.388889\\
0.0904724 -0.375\\
0.0871216 -0.361111\\
0.0837707 -0.347222\\
0.0804199 -0.3333335\\
0.0770691 -0.3194445\\
0.0737182 -0.3055555\\
0.0703674 -0.2916665\\
0.0670166 -0.277778\\
0.0636658 -0.263889\\
0.0603149 -0.25\\
0.0603149 -0.25\\
0.0536133 -0.222222\\
0.0469116 -0.1944445\\
0.04021 -0.1666665\\
0.0335083 -0.138889\\
0.0268066 -0.111111\\
0.020105 -0.0833335\\
0.0134033 -0.0555555\\
0.00670166 -0.0277778\\
0 0\\
-0.00670166 0.0277778\\
-0.0134033 0.0555555\\
-0.020105 0.0833335\\
-0.0268066 0.111111\\
-0.0335083 0.138889\\
-0.04021 0.1666665\\
-0.0469116 0.1944445\\
-0.0536133 0.222222\\
-0.0603149 0.25\\
-0.0603149 0.25\\
-0.0636658 0.263889\\
-0.0670166 0.277778\\
-0.0703674 0.2916665\\
-0.0737182 0.3055555\\
-0.0770691 0.3194445\\
-0.0804199 0.3333335\\
-0.0837707 0.347222\\
-0.0871216 0.361111\\
-0.0904724 0.375\\
-0.0938232 0.388889\\
-0.097174 0.402778\\
-0.100525 0.4166665\\
-0.103876 0.4305555\\
-0.107227 0.4444445\\
-0.110577 0.4583335\\
-0.113928 0.472222\\
-0.117279 0.486111\\
-0.12063 0.5\\
};
\addlegendentry{TSDT7};

\addplot [
color=black,
only marks,
mark=x,
mark options={solid}
]
table[row sep=crcr]{
0.120905 -0.5\\
0.117546 -0.486111\\
0.114188 -0.472222\\
0.110829 -0.4583335\\
0.107471 -0.4444445\\
0.104112 -0.4305555\\
0.100754 -0.4166665\\
0.0973955 -0.402778\\
0.094037 -0.388889\\
0.0906786 -0.375\\
0.0873201 -0.361111\\
0.0839616 -0.347222\\
0.0806032 -0.3333335\\
0.0772447 -0.3194445\\
0.0738862 -0.3055555\\
0.0705278 -0.2916665\\
0.0671693 -0.277778\\
0.0638108 -0.263889\\
0.0604524 -0.25\\
0.0604524 -0.25\\
0.0537354 -0.222222\\
0.0470185 -0.1944445\\
0.0403016 -0.1666665\\
0.0335847 -0.138889\\
0.0268677 -0.111111\\
0.0201508 -0.0833335\\
0.0134339 -0.0555555\\
0.00671693 -0.0277778\\
0 0\\
-0.00671693 0.0277778\\
-0.0134339 0.0555555\\
-0.0201508 0.0833335\\
-0.0268677 0.111111\\
-0.0335847 0.138889\\
-0.0403016 0.1666665\\
-0.0470185 0.1944445\\
-0.0537354 0.222222\\
-0.0604524 0.25\\
-0.0604524 0.25\\
-0.0638108 0.263889\\
-0.0671693 0.277778\\
-0.0705278 0.2916665\\
-0.0738862 0.3055555\\
-0.0772447 0.3194445\\
-0.0806032 0.3333335\\
-0.0839616 0.347222\\
-0.0873201 0.361111\\
-0.0906786 0.375\\
-0.094037 0.388889\\
-0.0973955 0.402778\\
-0.100754 0.4166665\\
-0.104112 0.4305555\\
-0.107471 0.4444445\\
-0.110829 0.4583335\\
-0.114188 0.472222\\
-0.117546 0.486111\\
-0.120905 0.5\\
};
\addlegendentry{FSDT5};

\end{axis}
\end{tikzpicture}%
}}
\subfigure[$\overline{w}(a/2,b/2,z)$]{\scalebox{0.75}{% This file was created by matlab2tikz v0.3.3.
% Copyright (c) 2008--2013, Nico Schlömer <nico.schloemer@gmail.com>
% All rights reserved.
% 
% The latest updates can be retrieved from
%   http://www.mathworks.com/matlabcentral/fileexchange/22022-matlab2tikz
% where you can also make suggestions and rate matlab2tikz.
% 
% 
% 
\begin{tikzpicture}

\begin{axis}[%
width=\figurewidth,
height=\figureheight,
scale only axis,
xmin=0.56,
xmax=0.76,
xlabel={$\overline{w}$},
ymin=-0.5,
ymax=0.5,
ylabel={Normalized thickness, z/h},
legend style={at={(0.35,0.5)},anchor=south west,draw=black,fill=white,legend cell align=left}
]
\addplot [
color=black,
only marks,
mark=o,
mark options={solid}
]
table[row sep=crcr]{
0.610616 -0.5\\
0.610004 -0.486111\\
0.609409 -0.472222\\
0.608832 -0.4583335\\
0.608272 -0.4444445\\
0.607729 -0.4305555\\
0.607203 -0.4166665\\
0.606694 -0.402778\\
0.606203 -0.388889\\
0.605729 -0.375\\
0.605272 -0.361111\\
0.604833 -0.347222\\
0.604411 -0.3333335\\
0.604006 -0.3194445\\
0.603618 -0.3055555\\
0.603247 -0.2916665\\
0.602894 -0.277778\\
0.602558 -0.263889\\
0.602239 -0.25\\
0.602239 -0.25\\
0.601653 -0.222222\\
0.601136 -0.1944445\\
0.600688 -0.1666665\\
0.600308 -0.138889\\
0.599998 -0.111111\\
0.599757 -0.0833335\\
0.599585 -0.0555555\\
0.599481 -0.0277778\\
0.599447 0\\
0.599481 0.0277778\\
0.599585 0.0555555\\
0.599757 0.0833335\\
0.599998 0.111111\\
0.600308 0.138889\\
0.600688 0.1666665\\
0.601136 0.1944445\\
0.601653 0.222222\\
0.602239 0.25\\
0.602239 0.25\\
0.602558 0.263889\\
0.602894 0.277778\\
0.603247 0.2916665\\
0.603618 0.3055555\\
0.604006 0.3194445\\
0.604411 0.3333335\\
0.604833 0.347222\\
0.605272 0.361111\\
0.605729 0.375\\
0.606203 0.388889\\
0.606694 0.402778\\
0.607203 0.4166665\\
0.607729 0.4305555\\
0.608272 0.4444445\\
0.608832 0.4583335\\
0.609409 0.472222\\
0.610004 0.486111\\
0.610616 0.5\\
};
\addlegendentry{HSDT13};

\addplot [
color=black,
dashed
]
table[row sep=crcr]{
0.61491 -0.5\\
0.61491 -0.486111\\
0.61491 -0.472222\\
0.61491 -0.4583335\\
0.61491 -0.4444445\\
0.61491 -0.4305555\\
0.61491 -0.4166665\\
0.61491 -0.402778\\
0.61491 -0.388889\\
0.61491 -0.375\\
0.61491 -0.361111\\
0.61491 -0.347222\\
0.61491 -0.3333335\\
0.61491 -0.3194445\\
0.61491 -0.3055555\\
0.61491 -0.2916665\\
0.61491 -0.277778\\
0.61491 -0.263889\\
0.61491 -0.25\\
0.61491 -0.25\\
0.61491 -0.222222\\
0.61491 -0.1944445\\
0.61491 -0.1666665\\
0.61491 -0.138889\\
0.61491 -0.111111\\
0.61491 -0.0833335\\
0.61491 -0.0555555\\
0.61491 -0.0277778\\
0.61491 0\\
0.61491 0.0277778\\
0.61491 0.0555555\\
0.61491 0.0833335\\
0.61491 0.111111\\
0.61491 0.138889\\
0.61491 0.1666665\\
0.61491 0.1944445\\
0.61491 0.222222\\
0.61491 0.25\\
0.61491 0.25\\
0.61491 0.263889\\
0.61491 0.277778\\
0.61491 0.2916665\\
0.61491 0.3055555\\
0.61491 0.3194445\\
0.61491 0.3333335\\
0.61491 0.347222\\
0.61491 0.361111\\
0.61491 0.375\\
0.61491 0.388889\\
0.61491 0.402778\\
0.61491 0.4166665\\
0.61491 0.4305555\\
0.61491 0.4444445\\
0.61491 0.4583335\\
0.61491 0.472222\\
0.61491 0.486111\\
0.61491 0.5\\
};
\addlegendentry{HSDT11A};

\addplot [
color=black,
only marks,
mark=diamond,
mark options={solid}
]
table[row sep=crcr]{
0.695349 -0.5\\
0.694875 -0.486111\\
0.694414 -0.472222\\
0.693966 -0.4583335\\
0.693532 -0.4444445\\
0.693111 -0.4305555\\
0.692704 -0.4166665\\
0.69231 -0.402778\\
0.691929 -0.388889\\
0.691561 -0.375\\
0.691207 -0.361111\\
0.690867 -0.347222\\
0.690539 -0.3333335\\
0.690225 -0.3194445\\
0.689925 -0.3055555\\
0.689637 -0.2916665\\
0.689363 -0.277778\\
0.689103 -0.263889\\
0.688856 -0.25\\
0.688856 -0.25\\
0.688401 -0.222222\\
0.688001 -0.1944445\\
0.687653 -0.1666665\\
0.687359 -0.138889\\
0.687119 -0.111111\\
0.686932 -0.0833335\\
0.686798 -0.0555555\\
0.686718 -0.0277778\\
0.686691 0\\
0.686718 0.0277778\\
0.686798 0.0555555\\
0.686932 0.0833335\\
0.687119 0.111111\\
0.687359 0.138889\\
0.687653 0.1666665\\
0.688001 0.1944445\\
0.688401 0.222222\\
0.688856 0.25\\
0.688856 0.25\\
0.689103 0.263889\\
0.689363 0.277778\\
0.689637 0.2916665\\
0.689925 0.3055555\\
0.690225 0.3194445\\
0.690539 0.3333335\\
0.690867 0.347222\\
0.691207 0.361111\\
0.691561 0.375\\
0.691929 0.388889\\
0.69231 0.402778\\
0.692704 0.4166665\\
0.693111 0.4305555\\
0.693532 0.4444445\\
0.693966 0.4583335\\
0.694414 0.472222\\
0.694875 0.486111\\
0.695349 0.5\\
};
\addlegendentry{HSDT11B};

\addplot [
color=black,
dash pattern=on 1pt off 3pt on 3pt off 3pt
]
table[row sep=crcr]{
0.696457 -0.5\\
0.696457 -0.486111\\
0.696457 -0.472222\\
0.696457 -0.4583335\\
0.696457 -0.4444445\\
0.696457 -0.4305555\\
0.696457 -0.4166665\\
0.696457 -0.402778\\
0.696457 -0.388889\\
0.696457 -0.375\\
0.696457 -0.361111\\
0.696457 -0.347222\\
0.696457 -0.3333335\\
0.696457 -0.3194445\\
0.696457 -0.3055555\\
0.696457 -0.2916665\\
0.696457 -0.277778\\
0.696457 -0.263889\\
0.696457 -0.25\\
0.696457 -0.25\\
0.696457 -0.222222\\
0.696457 -0.1944445\\
0.696457 -0.1666665\\
0.696457 -0.138889\\
0.696457 -0.111111\\
0.696457 -0.0833335\\
0.696457 -0.0555555\\
0.696457 -0.0277778\\
0.696457 0\\
0.696457 0.0277778\\
0.696457 0.0555555\\
0.696457 0.0833335\\
0.696457 0.111111\\
0.696457 0.138889\\
0.696457 0.1666665\\
0.696457 0.1944445\\
0.696457 0.222222\\
0.696457 0.25\\
0.696457 0.25\\
0.696457 0.263889\\
0.696457 0.277778\\
0.696457 0.2916665\\
0.696457 0.3055555\\
0.696457 0.3194445\\
0.696457 0.3333335\\
0.696457 0.347222\\
0.696457 0.361111\\
0.696457 0.375\\
0.696457 0.388889\\
0.696457 0.402778\\
0.696457 0.4166665\\
0.696457 0.4305555\\
0.696457 0.4444445\\
0.696457 0.4583335\\
0.696457 0.472222\\
0.696457 0.486111\\
0.696457 0.5\\
};
\addlegendentry{HSDT9};

\addplot [
color=black,
solid
]
table[row sep=crcr]{
0.741655 -0.5\\
0.741655 -0.486111\\
0.741655 -0.472222\\
0.741655 -0.4583335\\
0.741655 -0.4444445\\
0.741655 -0.4305555\\
0.741655 -0.4166665\\
0.741655 -0.402778\\
0.741655 -0.388889\\
0.741655 -0.375\\
0.741655 -0.361111\\
0.741655 -0.347222\\
0.741655 -0.3333335\\
0.741655 -0.3194445\\
0.741655 -0.3055555\\
0.741655 -0.2916665\\
0.741655 -0.277778\\
0.741655 -0.263889\\
0.741655 -0.25\\
0.741655 -0.25\\
0.741655 -0.222222\\
0.741655 -0.1944445\\
0.741655 -0.1666665\\
0.741655 -0.138889\\
0.741655 -0.111111\\
0.741655 -0.0833335\\
0.741655 -0.0555555\\
0.741655 -0.0277778\\
0.741655 0\\
0.741655 0.0277778\\
0.741655 0.0555555\\
0.741655 0.0833335\\
0.741655 0.111111\\
0.741655 0.138889\\
0.741655 0.1666665\\
0.741655 0.1944445\\
0.741655 0.222222\\
0.741655 0.25\\
0.741655 0.25\\
0.741655 0.263889\\
0.741655 0.277778\\
0.741655 0.2916665\\
0.741655 0.3055555\\
0.741655 0.3194445\\
0.741655 0.3333335\\
0.741655 0.347222\\
0.741655 0.361111\\
0.741655 0.375\\
0.741655 0.388889\\
0.741655 0.402778\\
0.741655 0.4166665\\
0.741655 0.4305555\\
0.741655 0.4444445\\
0.741655 0.4583335\\
0.741655 0.472222\\
0.741655 0.486111\\
0.741655 0.5\\
};
\addlegendentry{TSDT7};

\addplot [
color=black,
only marks,
mark=x,
mark options={solid}
]
table[row sep=crcr]{
0.738287 -0.5\\
0.738287 -0.486111\\
0.738287 -0.472222\\
0.738287 -0.4583335\\
0.738287 -0.4444445\\
0.738287 -0.4305555\\
0.738287 -0.4166665\\
0.738287 -0.402778\\
0.738287 -0.388889\\
0.738287 -0.375\\
0.738287 -0.361111\\
0.738287 -0.347222\\
0.738287 -0.3333335\\
0.738287 -0.3194445\\
0.738287 -0.3055555\\
0.738287 -0.2916665\\
0.738287 -0.277778\\
0.738287 -0.263889\\
0.738287 -0.25\\
0.738287 -0.25\\
0.738287 -0.222222\\
0.738287 -0.1944445\\
0.738287 -0.1666665\\
0.738287 -0.138889\\
0.738287 -0.111111\\
0.738287 -0.0833335\\
0.738287 -0.0555555\\
0.738287 -0.0277778\\
0.738287 0\\
0.738287 0.0277778\\
0.738287 0.0555555\\
0.738287 0.0833335\\
0.738287 0.111111\\
0.738287 0.138889\\
0.738287 0.1666665\\
0.738287 0.1944445\\
0.738287 0.222222\\
0.738287 0.25\\
0.738287 0.25\\
0.738287 0.263889\\
0.738287 0.277778\\
0.738287 0.2916665\\
0.738287 0.3055555\\
0.738287 0.3194445\\
0.738287 0.3333335\\
0.738287 0.347222\\
0.738287 0.361111\\
0.738287 0.375\\
0.738287 0.388889\\
0.738287 0.402778\\
0.738287 0.4166665\\
0.738287 0.4305555\\
0.738287 0.4444445\\
0.738287 0.4583335\\
0.738287 0.472222\\
0.738287 0.486111\\
0.738287 0.5\\
};
\addlegendentry{FSDT5};

\end{axis}
\end{tikzpicture}%
}}
\subfigure[$\overline{\sigma}_{xx}(a/2,b/2,z)$]{\scalebox{0.75}{% This file was created by matlab2tikz v0.3.3.
% Copyright (c) 2008--2013, Nico Schlömer <nico.schloemer@gmail.com>
% All rights reserved.
% 
% The latest updates can be retrieved from
%   http://www.mathworks.com/matlabcentral/fileexchange/22022-matlab2tikz
% where you can also make suggestions and rate matlab2tikz.
% 
% 
% 
\begin{tikzpicture}

\begin{axis}[%
width=\figurewidth,
height=\figureheight,
scale only axis,
xmin=-0.2,
xmax=0.15,
xlabel={$\overline{\sigma}_{xx}$},
ymin=-0.5,
ymax=0.5,
ylabel={Normalized thickness, z/h},
legend style={draw=black,fill=white,legend cell align=left},
legend pos = north west
]
\addplot [
color=black,
only marks,
mark=o,
mark options={solid}
]
table[row sep=crcr]{
-0.11331 -0.5\\
-0.103625 -0.486111\\
-0.094287 -0.472222\\
-0.0853269 -0.4583335\\
-0.0767523 -0.4444445\\
-0.0685712 -0.4305555\\
-0.0607912 -0.4166665\\
-0.0534196 -0.402778\\
-0.0464634 -0.388889\\
-0.0399294 -0.375\\
-0.0338243 -0.361111\\
-0.0281541 -0.347222\\
-0.022925 -0.3333335\\
-0.0181429 -0.3194445\\
-0.0138137 -0.3055555\\
-0.00994363 -0.2916665\\
-0.00654136 -0.277778\\
-0.00362976 -0.263889\\
-0.00264502 -0.25\\
-0.0305418 -0.25\\
-0.0272466 -0.222222\\
-0.0239166 -0.1944445\\
-0.0205563 -0.1666665\\
-0.01717 -0.138889\\
-0.013762 -0.111111\\
-0.0103367 -0.0833335\\
-0.00689833 -0.0555555\\
-0.00345133 -0.0277778\\
0 0\\
0.00345133 0.0277778\\
0.00689833 0.0555555\\
0.0103367 0.0833335\\
0.013762 0.111111\\
0.01717 0.138889\\
0.0205563 0.1666665\\
0.0239166 0.1944445\\
0.0272466 0.222222\\
0.0305418 0.25\\
0.00274919 0.25\\
0.00361392 0.263889\\
0.00652351 0.277778\\
0.00992408 0.2916665\\
0.0137925 0.3055555\\
0.0181202 0.3194445\\
0.0229007 0.3333335\\
0.0281282 0.347222\\
0.0337968 0.361111\\
0.0399004 0.375\\
0.0464328 0.388889\\
0.0533875 0.402778\\
0.0607576 0.4166665\\
0.0685361 0.4305555\\
0.0767158 0.4444445\\
0.0852889 0.4583335\\
0.0942476 0.472222\\
0.103584 0.486111\\
0.113289 0.5\\
};
\addlegendentry{HSDT13};

\addplot [
color=black,
dashed
]
table[row sep=crcr]{
-0.112079 -0.5\\
-0.102513 -0.486111\\
-0.0932861 -0.472222\\
-0.0844275 -0.4583335\\
-0.0759456 -0.4444445\\
-0.0678482 -0.4305555\\
-0.060143 -0.4166665\\
-0.0528374 -0.402778\\
-0.0459385 -0.388889\\
-0.039453 -0.375\\
-0.0333874 -0.361111\\
-0.0277481 -0.347222\\
-0.022541 -0.3333335\\
-0.0177718 -0.3194445\\
-0.0134459 -0.3055555\\
-0.00956864 -0.2916665\\
-0.00614527 -0.277778\\
-0.00318237 -0.263889\\
-0.00074845 -0.25\\
-0.021628 -0.25\\
-0.0193123 -0.222222\\
-0.0169657 -0.1944445\\
-0.0145922 -0.1666665\\
-0.0121955 -0.138889\\
-0.00977951 -0.111111\\
-0.00734813 -0.0833335\\
-0.00490518 -0.0555555\\
-0.00245452 -0.0277778\\
0 0\\
0.00245452 0.0277778\\
0.00490518 0.0555555\\
0.00734813 0.0833335\\
0.00977951 0.111111\\
0.0121955 0.138889\\
0.0145922 0.1666665\\
0.0169657 0.1944445\\
0.0193123 0.222222\\
0.021628 0.25\\
0.00074302 0.25\\
0.00316563 0.263889\\
0.00612689 0.277778\\
0.00954868 0.2916665\\
0.0134244 0.3055555\\
0.0177487 0.3194445\\
0.0225164 0.3333335\\
0.0277219 0.347222\\
0.0333597 0.361111\\
0.0394238 0.375\\
0.0459078 0.388889\\
0.0528052 0.402778\\
0.0601094 0.4166665\\
0.0678131 0.4305555\\
0.075909 0.4444445\\
0.0843895 0.4583335\\
0.0932467 0.472222\\
0.102472 0.486111\\
0.112058 0.5\\
};
\addlegendentry{HSDT11A};

\addplot [
color=black,
only marks,
mark=diamond,
mark options={solid}
]
table[row sep=crcr]{
-0.109336 -0.5\\
-0.0994172 -0.486111\\
-0.0900297 -0.472222\\
-0.0811748 -0.4583335\\
-0.072834 -0.4444445\\
-0.0649894 -0.4305555\\
-0.0576238 -0.4166665\\
-0.0507208 -0.402778\\
-0.0442644 -0.388889\\
-0.0382394 -0.375\\
-0.0326311 -0.361111\\
-0.0274257 -0.347222\\
-0.02261 -0.3333335\\
-0.0181714 -0.3194445\\
-0.0140986 -0.3055555\\
-0.0103813 -0.2916665\\
-0.0070129 -0.277778\\
-0.00400318 -0.263889\\
-0.00289686 -0.25\\
-0.034202 -0.25\\
-0.0301138 -0.222222\\
-0.0261273 -0.1944445\\
-0.0222297 -0.1666665\\
-0.0184083 -0.138889\\
-0.0146504 -0.111111\\
-0.0109433 -0.0833335\\
-0.00727439 -0.0555555\\
-0.00363084 -0.0277778\\
0 0\\
0.00363084 0.0277778\\
0.00727439 0.0555555\\
0.0109433 0.0833335\\
0.0146504 0.111111\\
0.0184083 0.138889\\
0.0222297 0.1666665\\
0.0261273 0.1944445\\
0.0301138 0.222222\\
0.034202 0.25\\
0.00300925 0.25\\
0.00398551 0.263889\\
0.0069937 0.277778\\
0.0103609 0.2916665\\
0.014077 0.3055555\\
0.0181487 0.3194445\\
0.022586 0.3333335\\
0.0274005 0.347222\\
0.0326046 0.361111\\
0.0382116 0.375\\
0.0442353 0.388889\\
0.0506904 0.402778\\
0.057592 0.4166665\\
0.0649562 0.4305555\\
0.0727993 0.4444445\\
0.0811387 0.4583335\\
0.0899921 0.472222\\
0.099378 0.486111\\
0.109316 0.5\\
};
\addlegendentry{HSDT11B};

\addplot [
color=black,
dash pattern=on 1pt off 3pt on 3pt off 3pt
]
table[row sep=crcr]{
-0.108408 -0.5\\
-0.0986209 -0.486111\\
-0.0893447 -0.472222\\
-0.0805824 -0.4583335\\
-0.0723172 -0.4444445\\
-0.064533 -0.4305555\\
-0.057214 -0.4166665\\
-0.0503451 -0.402778\\
-0.0439119 -0.388889\\
-0.0379004 -0.375\\
-0.0322974 -0.361111\\
-0.0270901 -0.347222\\
-0.0222662 -0.3333335\\
-0.0178144 -0.3194445\\
-0.0137235 -0.3055555\\
-0.00998334 -0.2916665\\
-0.00658459 -0.277778\\
-0.00352027 -0.263889\\
-0.0008547 -0.25\\
-0.0247565 -0.25\\
-0.0217957 -0.222222\\
-0.018909 -0.1944445\\
-0.0160873 -0.1666665\\
-0.0133211 -0.138889\\
-0.0106013 -0.111111\\
-0.00791853 -0.0833335\\
-0.00526357 -0.0555555\\
-0.00262715 -0.0277778\\
0 0\\
0.00262715 0.0277778\\
0.00526357 0.0555555\\
0.00791853 0.0833335\\
0.0106013 0.111111\\
0.0133211 0.138889\\
0.0160873 0.1666665\\
0.018909 0.1944445\\
0.0217957 0.222222\\
0.0247565 0.25\\
0.00084844 0.25\\
0.00350178 0.263889\\
0.00656495 0.277778\\
0.00996256 0.2916665\\
0.0137016 0.3055555\\
0.0177913 0.3194445\\
0.0222419 0.3333335\\
0.0270646 0.347222\\
0.0322707 0.361111\\
0.0378724 0.375\\
0.0438826 0.388889\\
0.0503145 0.402778\\
0.057182 0.4166665\\
0.0644996 0.4305555\\
0.0722825 0.4444445\\
0.0805462 0.4583335\\
0.089307 0.472222\\
0.0985817 0.486111\\
0.108388 0.5\\
};
\addlegendentry{HSDT9};

\addplot [
color=black,
solid
]
table[row sep=crcr]{
-0.100294 -0.5\\
-0.0921958 -0.486111\\
-0.084384 -0.472222\\
-0.0768767 -0.4583335\\
-0.0696742 -0.4444445\\
-0.0627763 -0.4305555\\
-0.056183 -0.4166665\\
-0.0498943 -0.402778\\
-0.0439103 -0.388889\\
-0.0382309 -0.375\\
-0.0328561 -0.361111\\
-0.0277859 -0.347222\\
-0.0230204 -0.3333335\\
-0.0185596 -0.3194445\\
-0.0144035 -0.3055555\\
-0.0105524 -0.2916665\\
-0.00700674 -0.277778\\
-0.0037691 -0.263889\\
-0.00091755 -0.25\\
-0.0266205 -0.25\\
-0.0236626 -0.222222\\
-0.0207048 -0.1944445\\
-0.017747 -0.1666665\\
-0.0147891 -0.138889\\
-0.0118313 -0.111111\\
-0.00887349 -0.0833335\\
-0.00591566 -0.0555555\\
-0.00295783 -0.0277778\\
0 0\\
0.00295783 0.0277778\\
0.00591566 0.0555555\\
0.00887349 0.0833335\\
0.0118313 0.111111\\
0.0147891 0.138889\\
0.017747 0.1666665\\
0.0207048 0.1944445\\
0.0236626 0.222222\\
0.0266205 0.25\\
0.00091078 0.25\\
0.0037493 0.263889\\
0.00698584 0.277778\\
0.0105304 0.2916665\\
0.0143805 0.3055555\\
0.0185355 0.3194445\\
0.0229953 0.3333335\\
0.0277598 0.347222\\
0.0328289 0.361111\\
0.0382026 0.375\\
0.043881 0.388889\\
0.049864 0.402778\\
0.0561516 0.4166665\\
0.0627438 0.4305555\\
0.0696407 0.4444445\\
0.0768422 0.4583335\\
0.0843484 0.472222\\
0.0921592 0.486111\\
0.100275 0.5\\
};
\addlegendentry{TSDT7};

\addplot [
color=black,
only marks,
mark=x,
mark options={solid}
]
table[row sep=crcr]{
-0.100515 -0.5\\
-0.0923987 -0.486111\\
-0.0845692 -0.472222\\
-0.0770451 -0.4583335\\
-0.0698264 -0.4444445\\
-0.062913 -0.4305555\\
-0.0563049 -0.4166665\\
-0.0500021 -0.402778\\
-0.0440047 -0.388889\\
-0.0383125 -0.375\\
-0.0329257 -0.361111\\
-0.0278442 -0.347222\\
-0.023068 -0.3333335\\
-0.0185972 -0.3194445\\
-0.0144319 -0.3055555\\
-0.0105722 -0.2916665\\
-0.00701865 -0.277778\\
-0.00377378 -0.263889\\
-0.00091506 -0.25\\
-0.0266013 -0.25\\
-0.0236456 -0.222222\\
-0.0206899 -0.1944445\\
-0.0177342 -0.1666665\\
-0.0147785 -0.138889\\
-0.0118228 -0.111111\\
-0.00886711 -0.0833335\\
-0.0059114 -0.0555555\\
-0.0029557 -0.0277778\\
0 0\\
0.0029557 0.0277778\\
0.0059114 0.0555555\\
0.00886711 0.0833335\\
0.0118228 0.111111\\
0.0147785 0.138889\\
0.0177342 0.1666665\\
0.0206899 0.1944445\\
0.0236456 0.222222\\
0.0266013 0.25\\
0.00090825 0.25\\
0.00375393 0.263889\\
0.0069977 0.277778\\
0.0105502 0.2916665\\
0.0144088 0.3055555\\
0.0185731 0.3194445\\
0.0230428 0.3333335\\
0.027818 0.347222\\
0.0328984 0.361111\\
0.0382842 0.375\\
0.0439753 0.388889\\
0.0499717 0.402778\\
0.0562734 0.4166665\\
0.0628805 0.4305555\\
0.0697928 0.4444445\\
0.0770105 0.4583335\\
0.0845336 0.472222\\
0.092362 0.486111\\
0.100496 0.5\\
};
\addlegendentry{FSDT5};

\end{axis}
\end{tikzpicture}%
}}
\subfigure[$\overline{\sigma}_{xz}(0,b/2,z)$]{\scalebox{0.75}{% This file was created by matlab2tikz v0.3.3.
% Copyright (c) 2008--2013, Nico Schlömer <nico.schloemer@gmail.com>
% All rights reserved.
% 
% The latest updates can be retrieved from
%   http://www.mathworks.com/matlabcentral/fileexchange/22022-matlab2tikz
% where you can also make suggestions and rate matlab2tikz.
% 
% 
% 
\begin{tikzpicture}

\begin{axis}[%
width=\figurewidth,
height=\figureheight,
scale only axis,
xmin=0,
xmax=0.012,
xlabel={$\overline{\sigma}_{xz}$},
ymin=-0.5,
ymax=0.5,
ylabel={Normalized thickness, z/h},
legend style={at={(0,0.4)},anchor=south west,draw=black,fill=white,legend cell align=left}
]
\addplot [
color=black,
only marks,
mark=o,
mark options={solid}
]
table[row sep=crcr]{
0 -0.5\\
0.00093251 -0.486111\\
0.00178347 -0.472222\\
0.00255608 -0.4583335\\
0.00325359 -0.4444445\\
0.00387934 -0.4305555\\
0.00443672 -0.4166665\\
0.00492918 -0.402778\\
0.00536025 -0.388889\\
0.0057335 -0.375\\
0.00605256 -0.361111\\
0.00632113 -0.347222\\
0.00654295 -0.3333335\\
0.0067218 -0.3194445\\
0.00686153 -0.3055555\\
0.00696603 -0.2916665\\
0.00703926 -0.277778\\
0.00708529 -0.263889\\
0.00710964 -0.25\\
0.00710964 -0.25\\
0.007796 -0.222222\\
0.00840385 -0.1944445\\
0.00893237 -0.1666665\\
0.00938082 -0.138889\\
0.00974858 -0.111111\\
0.0100351 -0.0833335\\
0.0102401 -0.0555555\\
0.0103632 -0.0277778\\
0.0104042 0\\
0.0103632 0.0277778\\
0.0102401 0.0555555\\
0.0100351 0.0833335\\
0.00974858 0.111111\\
0.00938082 0.138889\\
0.00893237 0.1666665\\
0.00840385 0.1944445\\
0.007796 0.222222\\
0.00710964 0.25\\
0.00710964 0.25\\
0.00708529 0.263889\\
0.00703926 0.277778\\
0.00696603 0.2916665\\
0.00686153 0.3055555\\
0.0067218 0.3194445\\
0.00654295 0.3333335\\
0.00632113 0.347222\\
0.00605256 0.361111\\
0.0057335 0.375\\
0.00536025 0.388889\\
0.00492918 0.402778\\
0.00443672 0.4166665\\
0.00387934 0.4305555\\
0.00325359 0.4444445\\
0.00255608 0.4583335\\
0.00178347 0.472222\\
0.00093251 0.486111\\
0 0.5\\
};
\addlegendentry{HSDT13};

\addplot [
color=black,
dashed
]
table[row sep=crcr]{
0 -0.5\\
0.00092245 -0.486111\\
0.0017643 -0.472222\\
0.00252866 -0.4583335\\
0.00321873 -0.4444445\\
0.00383776 -0.4305555\\
0.00438909 -0.4166665\\
0.0048761 -0.402778\\
0.00530224 -0.388889\\
0.00567101 -0.375\\
0.00598597 -0.361111\\
0.00625076 -0.347222\\
0.00646902 -0.3333335\\
0.00664449 -0.3194445\\
0.00678093 -0.3055555\\
0.00688215 -0.2916665\\
0.006952 -0.277778\\
0.00699438 -0.263889\\
0.00701335 -0.25\\
0.00701335 -0.25\\
0.00755625 -0.222222\\
0.00803735 -0.1944445\\
0.00845586 -0.1666665\\
0.00881111 -0.138889\\
0.00910255 -0.111111\\
0.00932971 -0.0833335\\
0.00949222 -0.0555555\\
0.00958983 -0.0277778\\
0.00962238 0\\
0.00958983 0.0277778\\
0.00949222 0.0555555\\
0.00932971 0.0833335\\
0.00910255 0.111111\\
0.00881111 0.138889\\
0.00845586 0.1666665\\
0.00803735 0.1944445\\
0.00755625 0.222222\\
0.00701335 0.25\\
0.00701335 0.25\\
0.00699438 0.263889\\
0.006952 0.277778\\
0.00688215 0.2916665\\
0.00678093 0.3055555\\
0.00664449 0.3194445\\
0.00646902 0.3333335\\
0.00625076 0.347222\\
0.00598597 0.361111\\
0.00567101 0.375\\
0.00530224 0.388889\\
0.0048761 0.402778\\
0.00438909 0.4166665\\
0.00383776 0.4305555\\
0.00321873 0.4444445\\
0.00252866 0.4583335\\
0.0017643 0.472222\\
0.00092245 0.486111\\
0 0.5\\
};
\addlegendentry{HSDT11A};

\addplot [
color=black,
only marks,
mark=diamond,
mark options={solid}
]
table[row sep=crcr]{
0 -0.5\\
0.00089823 -0.486111\\
0.00171359 -0.472222\\
0.00245073 -0.4583335\\
0.00311414 -0.4444445\\
0.00370816 -0.4305555\\
0.00423697 -0.4166665\\
0.0047046 -0.402778\\
0.00511495 -0.388889\\
0.00547179 -0.375\\
0.00577875 -0.361111\\
0.00603933 -0.347222\\
0.00625693 -0.3333335\\
0.00643483 -0.3194445\\
0.0065762 -0.3055555\\
0.00668414 -0.2916665\\
0.00676165 -0.277778\\
0.00681174 -0.263889\\
0.00683885 -0.25\\
0.00683885 -0.25\\
0.00761026 -0.222222\\
0.00828468 -0.1944445\\
0.00886443 -0.1666665\\
0.00935155 -0.138889\\
0.00974776 -0.111111\\
0.0100545 -0.0833335\\
0.0102728 -0.0555555\\
0.0104034 -0.0277778\\
0.0104469 0\\
0.0104034 0.0277778\\
0.0102728 0.0555555\\
0.0100545 0.0833335\\
0.00974776 0.111111\\
0.00935155 0.138889\\
0.00886443 0.1666665\\
0.00828468 0.1944445\\
0.00761026 0.222222\\
0.00683885 0.25\\
0.00683885 0.25\\
0.00681174 0.263889\\
0.00676165 0.277778\\
0.00668414 0.2916665\\
0.0065762 0.3055555\\
0.00643483 0.3194445\\
0.00625693 0.3333335\\
0.00603933 0.347222\\
0.00577875 0.361111\\
0.00547179 0.375\\
0.00511495 0.388889\\
0.0047046 0.402778\\
0.00423697 0.4166665\\
0.00370816 0.4305555\\
0.00311414 0.4444445\\
0.00245073 0.4583335\\
0.00171359 0.472222\\
0.00089823 0.486111\\
0 0.5\\
};
\addlegendentry{HSDT11B};

\addplot [
color=black,
dash pattern=on 1pt off 3pt on 3pt off 3pt
]
table[row sep=crcr]{
0 -0.5\\
0.00089085 -0.486111\\
0.00169981 -0.472222\\
0.00243138 -0.4583335\\
0.00308991 -0.4444445\\
0.00367961 -0.4305555\\
0.00420453 -0.4166665\\
0.00466862 -0.402778\\
0.00507566 -0.388889\\
0.00542935 -0.375\\
0.00573324 -0.361111\\
0.00599077 -0.347222\\
0.0062053 -0.3333335\\
0.00638004 -0.3194445\\
0.00651815 -0.3055555\\
0.00662266 -0.2916665\\
0.00669654 -0.277778\\
0.00674267 -0.263889\\
0.00676399 -0.25\\
0.00676399 -0.25\\
0.00738239 -0.222222\\
0.00792317 -0.1944445\\
0.00838817 -0.1666665\\
0.00877895 -0.138889\\
0.00909687 -0.111111\\
0.009343 -0.0833335\\
0.0095182 -0.0555555\\
0.00962308 -0.0277778\\
0.009658 0\\
0.00962308 0.0277778\\
0.0095182 0.0555555\\
0.009343 0.0833335\\
0.00909687 0.111111\\
0.00877895 0.138889\\
0.00838817 0.1666665\\
0.00792317 0.1944445\\
0.00738239 0.222222\\
0.00676399 0.25\\
0.00676399 0.25\\
0.00674267 0.263889\\
0.00669654 0.277778\\
0.00662266 0.2916665\\
0.00651815 0.3055555\\
0.00638004 0.3194445\\
0.0062053 0.3333335\\
0.00599077 0.347222\\
0.00573324 0.361111\\
0.00542935 0.375\\
0.00507566 0.388889\\
0.00466862 0.402778\\
0.00420453 0.4166665\\
0.00367961 0.4305555\\
0.00308991 0.4444445\\
0.00243138 0.4583335\\
0.00169981 0.472222\\
0.00089085 0.486111\\
0 0.5\\
};
\addlegendentry{HSDT9};

\addplot [
color=black,
solid
]
table[row sep=crcr]{
0 -0.5\\
0.00082808 -0.486111\\
0.00158788 -0.472222\\
0.00228204 -0.4583335\\
0.00291318 -0.4444445\\
0.00348389 -0.4305555\\
0.0039968 -0.4166665\\
0.00445453 -0.402778\\
0.00485967 -0.388889\\
0.00521484 -0.375\\
0.00552265 -0.361111\\
0.0057857 -0.347222\\
0.0060066 -0.3333335\\
0.00618796 -0.3194445\\
0.00633237 -0.3055555\\
0.00644244 -0.2916665\\
0.00652078 -0.277778\\
0.00657002 -0.263889\\
0.00659288 -0.25\\
0.00659288 -0.25\\
0.00726064 -0.222222\\
0.00784984 -0.1944445\\
0.00836047 -0.1666665\\
0.00879255 -0.138889\\
0.00914607 -0.111111\\
0.00942103 -0.0833335\\
0.00961743 -0.0555555\\
0.00973527 -0.0277778\\
0.00977455 0\\
0.00973527 0.0277778\\
0.00961743 0.0555555\\
0.00942103 0.0833335\\
0.00914607 0.111111\\
0.00879255 0.138889\\
0.00836047 0.1666665\\
0.00784984 0.1944445\\
0.00726064 0.222222\\
0.00659288 0.25\\
0.00659288 0.25\\
0.00657002 0.263889\\
0.00652078 0.277778\\
0.00644244 0.2916665\\
0.00633237 0.3055555\\
0.00618796 0.3194445\\
0.0060066 0.3333335\\
0.0057857 0.347222\\
0.00552265 0.361111\\
0.00521484 0.375\\
0.00485967 0.388889\\
0.00445453 0.402778\\
0.0039968 0.4166665\\
0.00348389 0.4305555\\
0.00291318 0.4444445\\
0.00228204 0.4583335\\
0.00158788 0.472222\\
0.00082808 0.486111\\
0 0.5\\
};
\addlegendentry{TSDT7};

\addplot [
color=black,
only marks,
mark=x,
mark options={solid}
]
table[row sep=crcr]{
0 -0.5\\
0.00082983 -0.486111\\
0.00159124 -0.472222\\
0.00228686 -0.4583335\\
0.00291931 -0.4444445\\
0.00349122 -0.4305555\\
0.00400519 -0.4166665\\
0.00446386 -0.402778\\
0.00486983 -0.388889\\
0.00522572 -0.375\\
0.00553414 -0.361111\\
0.00579771 -0.347222\\
0.00601903 -0.3333335\\
0.00620072 -0.3194445\\
0.00634539 -0.3055555\\
0.00645565 -0.2916665\\
0.00653411 -0.277778\\
0.00658339 -0.263889\\
0.00660625 -0.25\\
0.00660625 -0.25\\
0.00727255 -0.222222\\
0.00786046 -0.1944445\\
0.00836998 -0.1666665\\
0.00880111 -0.138889\\
0.00915386 -0.111111\\
0.00942822 -0.0833335\\
0.00962419 -0.0555555\\
0.00974177 -0.0277778\\
0.00978097 0\\
0.00974177 0.0277778\\
0.00962419 0.0555555\\
0.00942822 0.0833335\\
0.00915386 0.111111\\
0.00880111 0.138889\\
0.00836998 0.1666665\\
0.00786046 0.1944445\\
0.00727255 0.222222\\
0.00660625 0.25\\
0.00660625 0.25\\
0.00658339 0.263889\\
0.00653411 0.277778\\
0.00645565 0.2916665\\
0.00634539 0.3055555\\
0.00620072 0.3194445\\
0.00601903 0.3333335\\
0.00579771 0.347222\\
0.00553414 0.361111\\
0.00522572 0.375\\
0.00486983 0.388889\\
0.00446386 0.402778\\
0.00400519 0.4166665\\
0.00349122 0.4305555\\
0.00291931 0.4444445\\
0.00228686 0.4583335\\
0.00159124 0.472222\\
0.00082983 0.486111\\
0 0.5\\
};
\addlegendentry{FSDT5};

\end{axis}
\end{tikzpicture}%
}}
\caption{Displacements and stresses through the thickness for the square plates with simply supported edges with $a/h=$ 5, $h_H/h_f=$ 2 and $V_{\rm CN}^\ast=$ 0.28 and the CNT distribution is FG-X. The plate is subjected to a thermal load $T(x,y) = T_o\left( \frac{2z}{h} \right) \sin (\pi x/a) \sin (\pi y/b)$}
\label{fig:thickdispstresstherm}
\end{figure}

Through thickness variation of displacements and stresses for $h_H/h_f=$ 2 with $a/h =$ 5 and $a/b =$ 1 is demonstrated considering two locations in the sandwich plate in \frefs{fig:thickdispstressmech} and \ref{fig:thickdispstresstherm} for different structural models employed in the present investigation for mechanical and thermal load cases, respectively. A significant difference in the through thickness non-dimensionalized displacements and non-dimensionalized stresses can be seen while adopting different plate theories as noticed in \frefs{fig:thickdispstressmech} and \ref{fig:thickdispstresstherm}. The main difference among the higher-order theories is in accounting for the slope discontinuity in the in-plane displacements through the thickness whereas it is the inclusion of variation up to cubic term in the in-plane displacements among the lower-order theories. The influence of including the zig-zag function is clearly seen (HSDT13, HSDT11A) as there is a significant change in the material properties through the thickness. However, the shear stress variation predicted adopting HSDT11A and HSDT11B matches with some regions in the thickness direction in comparison with those of HSDT13.  It can be concluded that the performance of the structural model with the inclusion of zig-zag function is overall comparable with those of the higher-order model HSDT13 for the mechanical response of the sandwich structures. Similar experimentation of these types of theories for an applied thermal load is presented in \fref{fig:thickdispstresstherm}. The behaviour of HSDT11A and HSDT11B for thermal response is somewhat similar to those of mechanical case. However, it may be opined that the prediction of shear stresses by the structural model with the presence of linear/quadratic terms in the transverse response is overall very close to those of the higher-order model HSDT13 for the thermal response of the sandwich structures. 

%\begin{landscape}

\begin{table}[htpb]
\centering
\renewcommand\arraystretch{1.2}
\caption{Deflections and stresses for a simply supported square sandwich plates with homogeneous core and CNT reinforced face-sheets subjected to sinusoidally distributed load (mechanical load) with $a/h=$ 5 and 10. The volume fraction of the CNT, $V_{\rm CN}^\ast=$ 0.17 is assumed to functionally graded. The displacements and stresses are reported for the following locations: $\overline{u}=\overline{u}(0,b/2,h/2)$, $\overline{w} = \overline{w}(a/2,b/2,h/2)$, $\overline{\sigma}_{xx} = \overline{\sigma}_{xx}(a/2,a/2,h/2)$ and $\overline{\sigma}_{xz} = \overline{\sigma}_{xz}(0,b/2,0)$.}
\begin{tabular}{llrrrrrrrrr}
\hline
$a/h$ & Element & \multicolumn{4}{r}{$h_H/h_f=$ 2} && \multicolumn{4}{r}{$h_H/h_f=$ 6}\\
\cline{3-6}\cline{8-11}
& & $\overline{u}$ & $\overline{w}$ & $\overline{\sigma}_{xx}$  & $\overline{\sigma}_{xz}$ && $\overline{u}$ & $\overline{w}$ & $\overline{\sigma}_{xx}$  & $\overline{\sigma}_{xz}$ \\
\hline
\multirow{6}{*}{5}&	HSDT13	&	0.5463	&	0.1036	&	0.4476	&	0.4021	&&	0.5275	&	0.0533	&	0.4247	&	0.3154	\\
&	HSDT11A	&	0.5617	&	0.1084	&	0.4715	&	0.3985	&&	0.5390	&	0.0557	&	0.4491	&	0.3147	\\
&	HSDT11B	&	0.6701	&	0.0813	&	0.5496	&	0.3960	&&	0.6725	&	0.0504	&	0.5448	&	0.3168	\\
&	HSDT9	&	0.6883	&	0.0847	&	0.5747	&	0.3926	&&	0.6887	&	0.0526	&	0.5727	&	0.3162	\\
&	TSDT7	&	0.9108	&	0.0744	&	0.7578	&	0.3798	&&	0.6840	&	0.0523	&	0.5687	&	0.3155	\\
&	FSDT5	&	0.9042	&	0.0773	&	0.7525	&	0.3791	&&	0.6800	&	0.0538	&	0.5654	&	0.3149	\\
\cline{2-11}
\multirow{6}{*}{10} &	HSDT13	&	0.8132	&	0.0738	&	0.6717	&	0.3906	&&	0.6535	&	0.0470	&	0.5383	&	0.3183	\\
&	HSDT11A	&	0.8189	&	0.0749	&	0.6819	&	0.3891	&&	0.6574	&	0.0477	&	0.5466	&	0.3179	\\
&	HSDT11B	&	0.8663	&	0.0658	&	0.7155	&	0.3884	&&	0.6968	&	0.0461	&	0.5738	&	0.3188	\\
&	HSDT9	&	0.8721	&	0.0668	&	0.7255	&	0.3870	&&	0.7010	&	0.0468	&	0.5827	&	0.3184	\\
&	TSDT7	&	0.9359	&	0.0637	&	0.7779	&	0.3828	&&	0.6994	&	0.0467	&	0.5813	&	0.3181	\\
&	FSDT5	&	0.9341	&	0.0645	&	0.7765	&	0.3826	&&	0.6983	&	0.0471	&	0.5804	&	0.3179	\\
\hline
\end{tabular}
\label{table:vcnt17influence}
\end{table}
%\end{landscape}

%\begin{landscape}
\begin{table}[htpb]
\centering
\renewcommand\arraystretch{1.2}
\caption{Deflections and stresses for a simply supported square sandwich plates with homogeneous core and CNT reinforced face-sheets subjected to sinusoidally distributed load (mechanical load) with $a/h=$ 5 and 10. The volume fraction of the CNT, $V_{\rm CN}^\ast=$ 0.28 is assumed to functionally graded. The displacements and stresses are reported for the following locations: $\overline{u}=\overline{u}(0,b/2,h/2)$, $\overline{w} = \overline{w}(a/2,b/2,h/2)$, $\overline{\sigma}_{xx} = \overline{\sigma}_{xx}(a/2,a/2,h/2)$ and $\overline{\sigma}_{xz} = \overline{\sigma}_{xz}(0,b/2,0)$.}
\begin{tabular}{llrrrrrrrrr}
\hline
$a/h$ & Element & \multicolumn{4}{r}{$h_H/h_f=$ 2} && \multicolumn{4}{r}{$h_H/h_f=$ 6}\\
\cline{3-6}\cline{8-11}
& & $\overline{u}$ & $\overline{w}$ & $\overline{\sigma}_{xx}$  & $\overline{\sigma}_{xz}$ && $\overline{u}$ & $\overline{w}$ & $\overline{\sigma}_{xx}$ & $\overline{\sigma}_{xz}$ \\
\hline
\multirow{6}{*}{5} &	HSDT13	&	0.3603	&	0.0906	&	0.4809	&	0.3970	&&	0.4118	&	0.0481	&	0.5378	&	0.3167	\\
&	HSDT11A	&	0.3705	&	0.0945	&	0.5058	&	0.3933	&&	0.4201	&	0.0502	&	0.5691	&	0.3160	\\
&	HSDT11B	&	0.4343	&	0.0659	&	0.5800	&	0.3903	&&	0.5486	&	0.0441	&	0.7220	&	0.3186	\\
&	HSDT9	&	0.4456	&	0.0683	&	0.6051	&	0.3870	&&	0.5610	&	0.0458	&	0.7585	&	0.3180	\\
&	TSDT7	&	0.6314	&	0.0573	&	0.8542	&	0.3724	&&	0.5671	&	0.0455	&	0.7666	&	0.3167	\\
&	FSDT5	&	0.6260	&	0.0604	&	0.8472	&	0.3715	&&	0.5626	&	0.0471	&	0.7606	&	0.3158	\\
\cline{2-11}
\multirow{6}{*}{10}&	HSDT13	&	0.5584	&	0.0574	&	0.7496	&	0.3846	&&	0.5375	&	0.0404	&	0.7190	&	0.3205	\\
&	HSDT11A	&	0.5624	&	0.0583	&	0.7614	&	0.3830	&&	0.5406	&	0.0410	&	0.7307	&	0.3201	\\
&	HSDT11B	&	0.5940	&	0.0485	&	0.7973	&	0.3821	&&	0.5801	&	0.0391	&	0.7764	&	0.3212	\\
&	HSDT9	&	0.5979	&	0.0492	&	0.8086	&	0.3807	&&	0.5835	&	0.0397	&	0.7885	&	0.3208	\\
&	TSDT7	&	0.6519	&	0.0458	&	0.8810	&	0.3760	&&	0.5848	&	0.0396	&	0.7902	&	0.3202	\\
&	FSDT5	&	0.6504	&	0.0466	&	0.8791	&	0.3757	&&	0.5835	&	0.0400	&	0.7885	&	0.3199	\\
\hline
\end{tabular}
\label{table:vcnt28influence}
\end{table}

%\end{landscape}

% figures...FGM 121 displacement and stress plot through the thickness

% thermal loading results
% Thermal loading

%\begin{landscape}

\begin{table}[htpb]
\centering
\renewcommand\arraystretch{1.2}
\caption{Deflections and stresses for a simply supported square sandwich plates with homogeneous core and CNT reinforced face-sheets subjected to sinusoidally distributed load (thermal load) with $a/h=$ 5 and 10. The volume fraction of the CNT, $V_{\rm CN}^\ast=$ 0.17 is assumed to functionally graded. The displacements and stresses are reported for the following locations: $\overline{u}=\overline{u}(0,b/2,h/2)$, $\overline{w} = \overline{w}(a/2,b/2,h/2)$, $\overline{\sigma}_{xx} = \overline{\sigma}_{xx}(a/2,a/2,h/2)$ and $\overline{\sigma}_{xz} = \overline{\sigma}_{xz}(0,b/2,-h/6)$.}
\begin{tabular}{llrrrrrrrrr}
\hline
$a/h$ & Element & \multicolumn{4}{r}{$h_H/h_f=$ 2} && \multicolumn{4}{r}{$h_H/h_f=$ 6}\\
\cline{3-6}\cline{8-11}
& & $\overline{u}$ & $\overline{w}$ & $\overline{\sigma}_{xx}$ & $\overline{\sigma}_{xz}$ && $\overline{u}$ & $\overline{w}$ & $\overline{\sigma}_{xx}$ & $\overline{\sigma}_{xz}$ \\
\hline
\multirow{6}{*}{5}&	HSDT13	&	0.1357	&	0.6106	&	0.1133	&	0.0893	&&	0.0907	&	0.3786	&	0.0757	&	0.0620	\\
&	HSDT11A	&	0.1350	&	0.6149	&	0.1121	&	0.0846	&&	0.0901	&	0.3755	&	0.0748	&	0.0532	\\
&	HSDT11B	&	0.1308	&	0.6953	&	0.1093	&	0.0886	&&	0.0730	&	0.4137	&	0.0610	&	0.0612	\\
&	HSDT9	&	0.1305	&	0.6965	&	0.1084	&	0.0839	&&	0.0730	&	0.4108	&	0.0607	&	0.0525	\\
&	TSDT7	&	0.1206	&	0.7417	&	0.1003	&	0.0836	&&	0.0682	&	0.4089	&	0.0567	&	0.0525	\\
&	FSDT5	&	0.1209	&	0.7383	&	0.1005	&	0.0837	&&	0.0688	&	0.4075	&	0.0571	&	0.0527	\\
\cline{2-11}
\multirow{6}{*}{10}&	HSDT13	&	0.1242	&	0.7070	&	0.1038	&	0.0308	&&	0.0725	&	0.4035	&	0.0606	&	0.0159	\\
&	HSDT11A	&	0.1243	&	0.7122	&	0.1033	&	0.0306	&&	0.0725	&	0.4042	&	0.0603	&	0.0147	\\
&	HSDT11B	&	0.1222	&	0.7376	&	0.1021	&	0.0303	&&	0.0672	&	0.4143	&	0.0562	&	0.0155	\\
&	HSDT9	&	0.1224	&	0.7409	&	0.1018	&	0.0301	&&	0.0673	&	0.4149	&	0.0560	&	0.0143	\\
&	TSDT7	&	0.1196	&	0.7545	&	0.0994	&	0.0298	&&	0.0661	&	0.4144	&	0.0550	&	0.0141	\\
&	FSDT5	&	0.1197	&	0.7536	&	0.0995	&	0.0298	&&	0.0663	&	0.4140	&	0.0551	&	0.0142	\\
\hline
\end{tabular}
\label{table:vcnt17influencetherm}
\end{table}
%\end{landscape}

%\begin{landscape}
\begin{table}[htpb]
\centering
\renewcommand\arraystretch{1.2}
\caption{Deflections and stresses for a simply supported square sandwich plates with homogeneous core and CNT reinforced face-sheets subjected to sinusoidally distributed load (thermal load) with $a/h=$ 5 and 10. The volume fraction of the CNT, $V_{\rm CN}^\ast=$ 0.28 is assumed to functionally graded. The displacements and stresses are reported for the following locations: $\overline{u}=\overline{u}(0,b/2,h/2)$, $\overline{w} = \overline{w}(a/2,b/2,h/2)$, $\overline{\sigma}_{xx} = \overline{\sigma}_{xx}(a/2,a/2,h/2)$ and $\overline{\sigma}_{xz} = \overline{\sigma}_{xz}(0,b/2,-h/6)$.}
\begin{tabular}{llrrrrrrrrr}
\hline
$a/h$ & Element & \multicolumn{4}{r}{$h_H/hf=$ 2} && \multicolumn{4}{r}{$h_H/hf=$ 6}\\
\cline{3-6}\cline{8-11}
& & $\overline{u}$ & $\overline{w}$ & $\overline{\sigma}_{xx}$ & $\overline{\sigma}_{xz}$ && $\overline{u}$ & $\overline{w}$ & $\overline{\sigma}_{xx}$ & $\overline{\sigma}_{xz}$ \\
\hline
\multirow{6}{*}{5}&	HSDT13	&	0.1145	&	0.5642	&	0.1552	&	0.1155	&&	0.0870	&	0.3884	&	0.1179	&	0.0768	\\
&	HSDT11A	&	0.1142	&	0.5672	&	0.1542	&	0.1111	&&	0.0868	&	0.3863	&	0.1172	&	0.0683	\\
&	HSDT11B	&	0.1123	&	0.6351	&	0.1523	&	0.1149	&&	0.0748	&	0.4246	&	0.1014	&	0.0761	\\
&	HSDT9	&	0.1121	&	0.6348	&	0.1515	&	0.1105	&&	0.0749	&	0.4228	&	0.1012	&	0.0675	\\
&	TSDT7	&	0.1076	&	0.6615	&	0.1454	&	0.1104	&&	0.0709	&	0.4237	&	0.0958	&	0.0675	\\
&	FSDT5	&	0.1077	&	0.6581	&	0.1456	&	0.1105	&&	0.0713	&	0.4217	&	0.0964	&	0.0677	\\
\cline{2-11}
\multirow{6}{*}{10}&	HSDT13	&	0.1090	&	0.6379	&	0.1478	&	0.0438	&&	0.0736	&	0.4186	&	0.0998	&	0.0229	\\
&	HSDT11A	&	0.1092	&	0.6425	&	0.1476	&	0.0436	&&	0.0738	&	0.4199	&	0.0997	&	0.0217	\\
&	HSDT11B	&	0.1079	&	0.6633	&	0.1464	&	0.0433	&&	0.0698	&	0.4304	&	0.0946	&	0.0224	\\
&	HSDT9	&	0.1082	&	0.6662	&	0.1463	&	0.0432	&&	0.0700	&	0.4315	&	0.0946	&	0.0212	\\
&	TSDT7	&	0.1069	&	0.6744	&	0.1445	&	0.0430	&&	0.0690	&	0.4318	&	0.0932	&	0.0211	\\
&	FSDT5	&	0.1069	&	0.6735	&	0.1446	&	0.0430	&&	0.0691	&	0.4312	&	0.0934	&	0.0211	\\
\hline
\end{tabular}
\label{table:vcnt28influencetherm}
\end{table}

%\end{landscape}

% figures...FGM 121 displacement and stress plot through the thickness

%\begin{table}
%\label{table:cntproperty}
%\end{table}

\subsection{Free flexural vibrations}
Next, the free vibration characteristics of sandwich plate is numerically studied. In all cases, we present the non-dimensionalized free flexural frequency defined as:
\begin{equation}
\Omega = \frac{\omega a^2}{h} \sqrt{ \frac{ \rho_{_H}}{E_{_H}}}
\label{eqn:vibnondim}
\end{equation}
where $\omega$ is the natural frequency, $\rho_{_H}, E_{_H}$ are the mass density and the Young's modulus of the homogeneous core evaluated at $T_o=$ 300K. Before proceeding with a detailed numerical study on the effect of CNT volume fraction, temperature effects, the core thickness and the plate thickness on the fundamental frequency, the formulation developed herein is validated against available results in the literature~\cite{wangshen2012}, which is based on the first order shear deformation theory. Table \ref{vibvalid} presents the non-dimensionalized frequency of sandwich plates with CNT reinforced facesheets in thermal environment based on various plate theories. For this study, a 8$\times$8 structures Q8 mesh was found to be adequate to model the entire plate. It is seen that the results from the present formulation are in very good agreement with the results available in the literature~\cite{wangshen2012}. The effect of volume fraction of the CNTs, $V_{\rm CN}^\ast$ and core-to-facesheet thickness ratio is also shown. With increasing temperature, the non-dimensionalized frequency decreases, whilst it increases with increasing volume fraction and core-to-facesheet thickness as expected.

% convergence of Q8-HSDT13 element with mesh 
\begin{table}[htpb]
\centering
\renewcommand\arraystretch{1.2}
\caption{Comparison of fundamental frequency $\Omega = \omega\frac{a^2}{h}\sqrt{ \frac{ \rho_{_H}}{E_{_H}}}$ for sandwich plates with CNT reinforced facesheets in thermal environment with $a/b=$ 1 and $a/h=$ 20 for various core-to-facesheet thickness $(h_H/h_f=$4,6$)$ and CNT volume fraction $V_{\rm CNT}^\ast$. The CNTs are assumed to the functionally graded, i.e., FGX.}
\begin{tabular}{clrrrrrrr}
\hline
Temperature & Theory & \multicolumn{3}{r}{$h_H/h_f=$ 4} && \multicolumn{3}{r}{$h_H/h_f=$ 6} \\
\cline{3-5}\cline{7-9}
& & $V_{\rm CN}^\ast=$ 0.12 & $V_{\rm CN}^\ast=$ 0.17 & $V_{\rm CN}^\ast=$ 0.28 &&  $V_{\rm CN}^\ast=$ 0.12 & $V_{\rm CN}^\ast=$ 0.17 & $V_{\rm CN}^\ast=$ 0.28 \\
\hline
\multirow{7}{*}{300}&	HSDT13	&	4.6518	&	5.0381	&	5.6422	&&	4.8992	&	5.1753	&	5.6278	\\
&	HSDT11A	&	4.6508	&	5.0366	&	5.6394	&&	4.8980	&	5.1737	&	5.6253	\\
&	HSDT11B	&	4.6803	&	5.0683	&	5.6980	&&	4.9111	&	5.1886	&	5.6534	\\
&	HSDT9	&	4.6793	&	5.0669	&	5.6955	&&	4.9099	&	5.1871	&	5.6510	\\
&	TSDT7	&	4.6808	&	5.0697	&	5.7025	&&	4.9111	&	5.1881	&	5.6524	\\
&	FSDT5	&	4.6765	&	5.0639	&	5.6935	&&	4.9066	&	5.1826	&	5.6448	\\
&	Ref.~\cite{wangshen2012} &	4.6845	&	5.0763	&	5.7131	&&	4.9119	&	5.1905	&	5.6569	\\
\cline{2-9}
\multirow{7}{*}{500}&	HSDT13	&	4.4425	&	4.8244	&	5.4238	&&	4.6662	&	4.9408	&	5.3956	\\
&	HSDT11A	&	4.4416	&	4.8232	&	5.4216	&&	4.6651	&	4.9395	&	5.3936	\\
&	HSDT11B	&	4.4809	&	4.8658	&	5.5005	&&	4.6823	&	4.9590	&	5.4308	\\
&	HSDT9	&	4.4801	&	4.8647	&	5.4985	&&	4.6813	&	4.9577	&	5.4289	\\
&	TSDT7	&	4.4819	&	4.8681	&	5.5075	&&	4.6824	&	4.9587	&	5.4305	\\
&	FSDT5	&	4.4776	&	4.8623	&	5.4982	&&	4.6779	&	4.9533	&	5.4227	\\
&	Ref.~\cite{wangshen2012}	&	4.4853	&	4.8743	&	5.5180	&&	4.6831	&	4.9609	&	5.4350	\\
\hline
\end{tabular}
\label{vibvalid}
\end{table}

% effect of a/h, volume fraction and type of reinforcement
\begin{table}[htpb]
\centering
\renewcommand\arraystretch{1.2}
\caption{Influence of the plate thickness $a/h$, the volume fraction of the nanotubes $V_{\rm CN}^ast$ and core-to-the thickness ratio $h_H/h_f$ on the fundamental frequency parameter $\Omega$ for a simply supported square sandwich plate with a homogeneous core. The temperature is assumed to be $T_o=$ 300K. The CNTs are assumed to the functionally graded, i.e., FGX.}
\begin{tabular}{clrrrrrrr}
\hline
$a/h$ & Element & \multicolumn{3}{r}{$h_H/h_f=$ 2} && \multicolumn{3}{r}{$h_H/h_f=$ 6}\\
\cline{3-5} \cline{7-9}
& Type & $V_{\rm CN}^\ast=$ 0.12 &  $V_{\rm CN}^\ast=$ 0.17 & $V_{\rm CN}^\ast=$ 0.28 && $V_{\rm CN}^\ast=$ 0.12 & $V_{\rm CN}^\ast=$ 0.17 & $V_{\rm CN}^\ast=$ 0.28 \\
\hline
\multirow{6}{*}{5} &	HSDT13	&	3.3961	&	3.7815	&	4.0457	&&	4.4193	&	4.6221	&	4.8681	\\
&	HSDT11A	&	3.3858	&	3.7714	&	4.0327	&&	4.4058	&	4.6097	&	4.8560	\\
&	HSDT11B	&	3.8434	&	4.2602	&	4.7377	&&	4.5334	&	4.7512	&	5.0871	\\
&	HSDT9	&	3.8381	&	4.2554	&	4.7318	&&	4.5212	&	4.7402	&	5.0774	\\
&	TSDT7	&	4.0745	&	4.5394	&	5.1622	&&	4.5365	&	4.7534	&	5.0935	\\
&	FSDT5	&	4.0164	&	4.4561	&	5.0308	&&	4.4820	&	4.6891	&	5.0097	\\
\cline{2-9}
\multirow{6}{*}{10}& HSDT13	&	4.0191	&	4.5719	&	5.1744	&&	4.7838	&	5.0419	&	5.4336	\\
&	HSDT11A	&	4.0166	&	4.5688	&	5.1689	&&	4.7798	&	5.0378	&	5.4289	\\
&	HSDT11B	&	4.2753	&	4.8382	&	5.6273	&&	4.8264	&	5.0895	&	5.5223	\\
&	HSDT9	&	4.2739	&	4.8364	&	5.6240	&&	4.8225	&	5.0857	&	5.5182	\\
&	TSDT7	&	4.3659	&	4.9522	&	5.8245	&&	4.8271	&	5.0897	&	5.5232	\\
&	FSDT5	&	4.3472	&	4.9238	&	5.7744	&&	4.8100	&	5.0691	&	5.4952	\\
\hline
\end{tabular}
\label{effectofah}
\end{table}

% effect of temperature, type reinforcement and volume fraction
\begin{table}[htpb]
\centering
\renewcommand\arraystretch{1.2}
\caption{Influence of the temperature $T$, plate thickness $a/h$ and the type of CNT distribution, viz., UD/FGX on the fundamental frequency parameter $\Omega$ for a simply supported square plate with a homogeneous core with $h_H/h_f =$ 2 and $V_{\rm CN}^\ast=$ 0.17.}
\begin{tabular}{clrrrrr}
\hline
$a/h$ & Element & \multicolumn{2}{l}{UD} && \multicolumn{2}{l}{FGX}\\
\cline{3-4} \cline{6-7}
& Type & $T=$ 300K & $T=$ 500K && $T=$300K & $T=$ 500K\\
\hline
\multirow{6}{*}{5}&	HSDT13	&	3.8203	&	3.5588	&&	3.7815	&	3.4910	\\
&	HSDT11A	&	3.8108	&	3.5492	&&	3.7714	&	3.4809	\\
&	HSDT11B	&	4.1736	&	3.9680	&&	4.2602	&	4.0277	\\
&	HSDT9	&	4.1677	&	3.9623	&&	4.2554	&	4.0232	\\
&	TSDT7	&	4.3199	&	4.1501	&&	4.5394	&	4.3657	\\
&	FSDT5	&	4.2504	&	4.0789	&&	4.4561	&	4.2801	\\
\cline{2-7}
\multirow{6}{*}{10}&	HSDT13	&	4.4397	&	4.2129	&&	4.5719	&	4.3092	\\
&	HSDT11A	&	4.4373	&	4.2105	&&	4.5688	&	4.3063	\\
&	HSDT11B	&	4.6095	&	4.4305	&&	4.8382	&	4.6441	\\
&	HSDT9	&	4.6078	&	4.4291	&&	4.8364	&	4.6427	\\
&	TSDT7	&	4.6655	&	4.5034	&&	4.9522	&	4.7899	\\
&	FSDT5	&	4.6426	&	4.4795	&&	4.9238	&	4.7602	\\
\hline
\end{tabular}
\label{effectoftemperature}
\end{table}

% include figures here
\begin{figure}[htpb]
\centering
\includegraphics[scale=0.6]{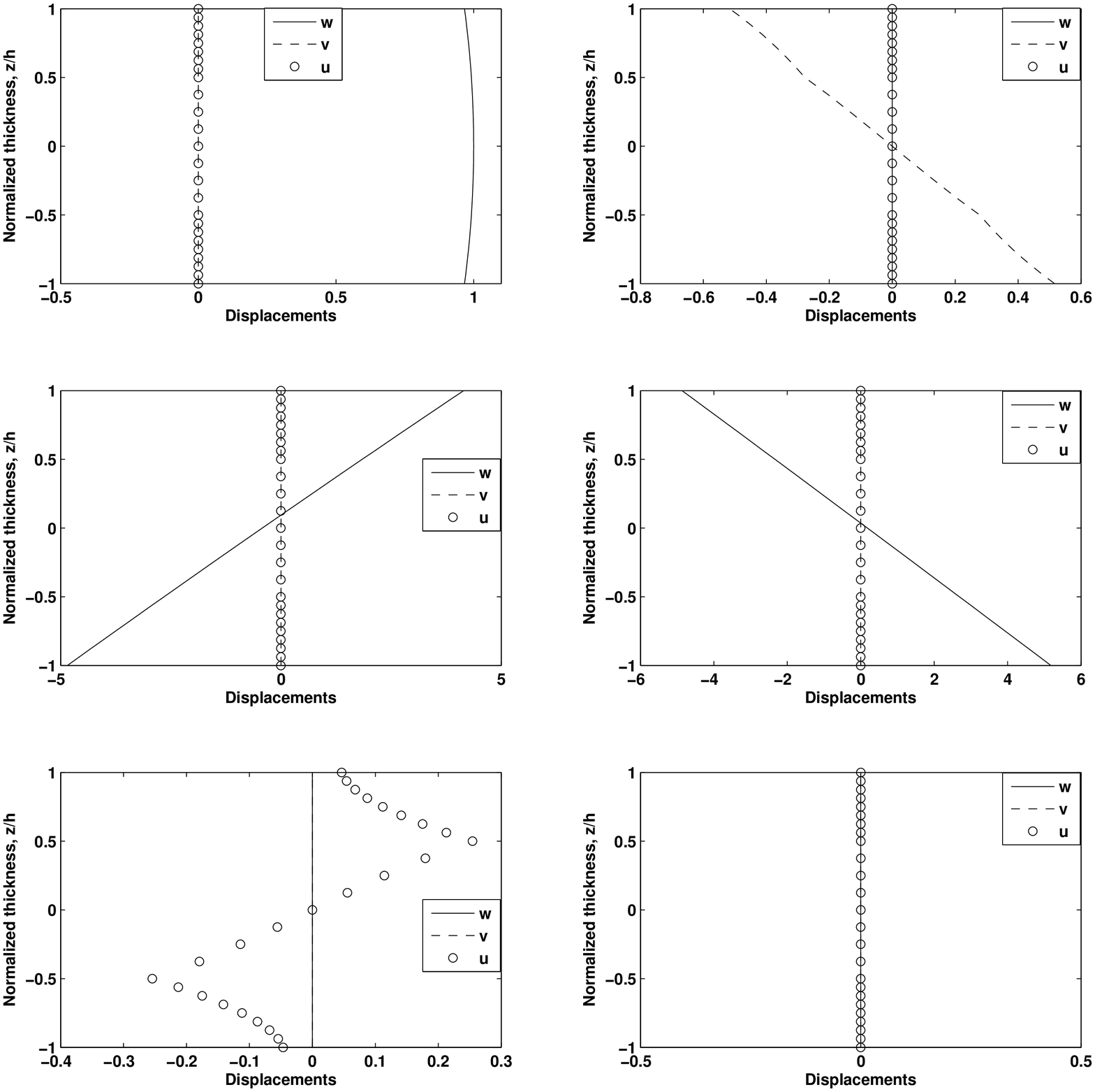}
\caption{Deflected shapes $u_i(a/2,a/2,z),~i=x,y,z$ of the six modes for the square plates with simply supported edges with thickness, $a/h=$ 5.}
\label{fig:centerA}
\end{figure}

% figures...FGM 121 n = 2 at the quarter of the plate
\begin{figure}[htpb]
\centering
\includegraphics[scale=0.6]{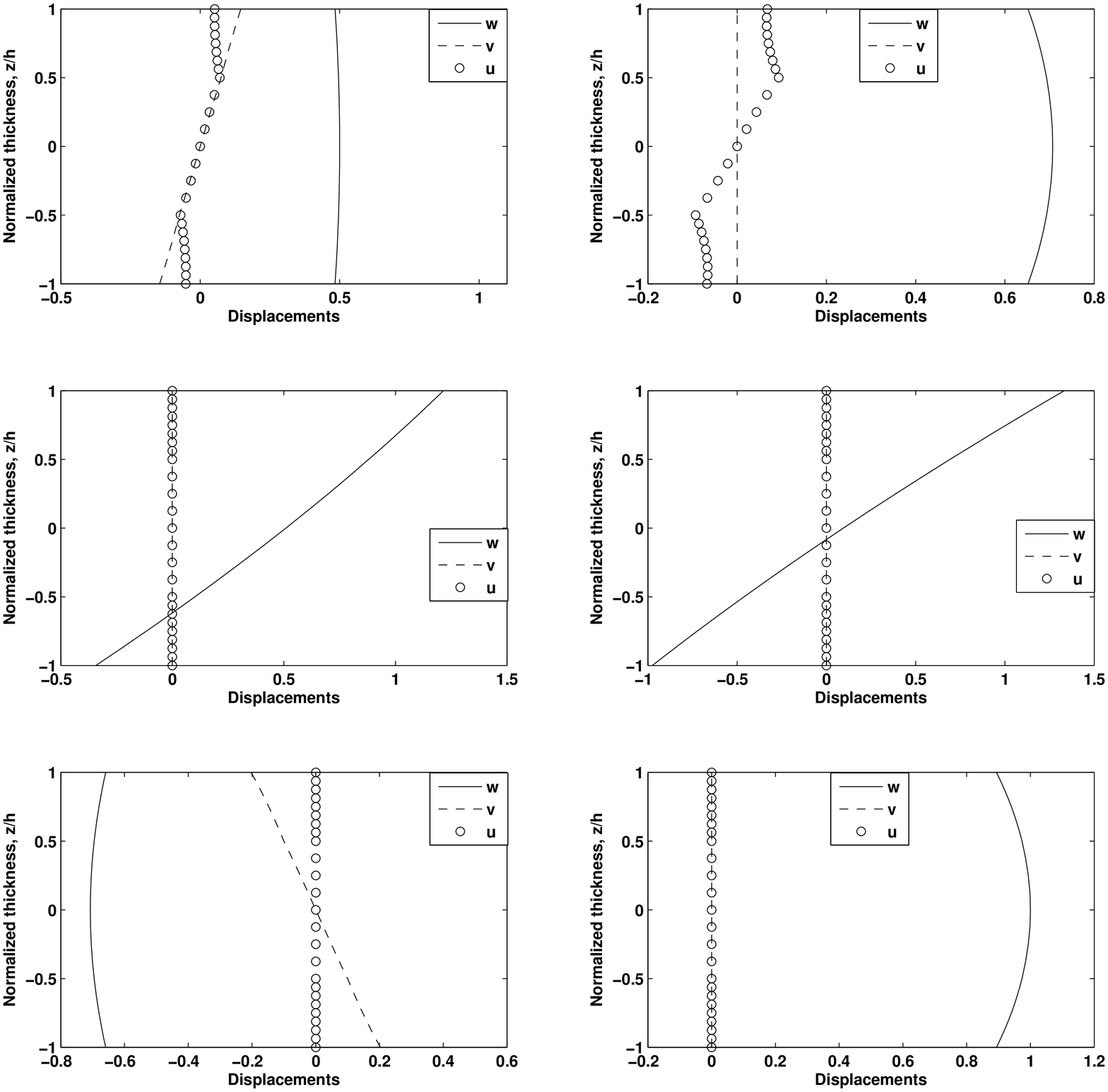}
\caption{Deflected shapes $u_i(3a/4,3a/4,z),~i=x,y,z$ of the six modes for the square plates with simply supported edges with thickness, $a/h=$ 5.}
\label{fig:quatA}
\end{figure}

Next, the influence of the plate side-to-thickness ratio, $a/h$, the volume fraction of the CNTs $V_{\rm CN}^\ast$, the CNT distribution, the core-to-facesheet thickness ratio $h_H/h_f$ and the temperature on the fundamental frequency is numerically studied by employing various plate theories. The results are tabulated in Table \ref{effectofah}. It can be opined from Table \ref{effectofah} that the performance of higher-order theories, HSDT13 and HSDT11A are comparable as demonstrated in static bending case. The CNTs are reinforced functionally graded through the thickness according to \Erefs{eqn:topsh} - (\ref{eqn:botsh}). It can be seen that the non-dimensionalized frequency $(\Omega)$ of the plate increases with decrease in the plate thickness, whilst increases with increasing CNT volume fraction and the core-to-facesheet thickness. This can be attributed to the local flexibility due to decreasing plate thickness and increasing the homogeneous core thickness. The influence of the temperature and the type of CNT distribution, viz, UD/FGX on the non-dimensionalized frequency for a square plate with $a/h=$ 5 is shown in Table \ref{effectoftemperature}. For all the plate theories considered in this study, it is seen that the non-dimenisonalized frequency increases with the increase in the homogeneous core thickness and decreases with increasing temperature. It can be observed from Table \ref{effectoftemperature} that for $a/h=$ 10, the non-dimesninoalized frequency increases when the CNT distribution is changed from UD to FGX, irrespective of the structural model employed. However, for $a/h=$ 5, structural models HSDT11B, HSDT9, TSDT7, FSDT5 predicts increasing fundamental frequency when changing the CNT distribution from UD to FGX, whilst HSDT13 and HSDT11A predicts decrease in the frequency.

\frefs{fig:centerA} - \ref{fig:quatA} show the relative displacements through the thickness of the sandwich plate for the first six modes using HSDT13. The plate is square simply supported with $h_H/h_f=$ 2 and $a/h =$ 5. The displacements $(u, v, w)$ are plotted along the lines $(a/2, b/2, z)$ and $(3a/4, 3b/4, z)$, where $−h/2 \le  z \le h/2$. It can be seen from the flexural modes in \frefs{fig:centerA} - \ref{fig:quatA} where the transverse displacement $w$ is not uniform at the chosen locations exhibit the existence of normal stresses in the thickness direction. In flexural modes 2 -- 5 shown in \fref{fig:centerA}, the deflected shape retains the thickness at $(a/2, b/2, z)$ whereas the thickness of the plate is compressed at the other location, i.e., at $(3a/4, 3b/4, z)$ as seen by the modes 1, 2 \& 5 in \fref{fig:quatA}. It can be also viewed that the variation of in-plane displacement, in general, is linear or nonlinear, irrespective of the types of modes. Furthermore, there may be a possibility of occurring of extensional mode at some locations. In general, it may be opined that the types of modes in the thickness direction depend on the location in the given sandwich structures.

\section{Conclusions}
\label{conclu}
The static and dynamic responses of sandwich plates with CNT reinforced facesheets are studied considering various parameters such as the sandwich type, temperature effects, the thickness ratio and the volume fraction of the CNT. Different plate models are employed in predicting the global structural behaviour and their through thickness variations in the sandwich plate. From a detailed parametric investigation on the effectiveness of higher-order models (HSDT13, HSDT11A and HSDT11B) over the other theories, neglecting the stretching terms in transverse displacement function, the following observations can be made:
\begin{itemize}
\item HSDT11A and HSDT11B, to some extent predict the global response of stuctures, irrespective of the type of load, viz., mechanical or thermal. However, there is a variation in evaluating the shear stress through the thickness of the sandwich plates.
\item HSDT13 is capable of evaluating the displacements and stresses accurately for the global as well as the local response characteristics of the structure.
\item The performance of the higher-order model, HSDT13 for the thick plate is significantly different over other lower theories considering in the present analysis.
\item In-plane stress variation is nonlinear and has discontinuity at the layer interface.
\item Increase in the volume fraction of CNT distribution in the facesheet, in general, decreases the deflection.
\item The effect of temperature variation affects the stress distribution noticeable in comparison with those of the mechanical load.
\item The variation in numerically computed natural frequencies through different plate theories is significant.
\item For thick plate, change in the CNT distribution from UD to FGX, the higher order structural models, viz., HSDT13 and HSDT11A predicts decrease in the non-dimensionalized frequency compared to other structural models.
\item The occurrence of type of flexural/extensional modes in the thickness direction depends on the location of the structures.
\item The exhibition of normal stress in the thickness direction through mode shape can be demonstrated.
\end{itemize}

\section*{Acknowledgement}
S Natarajan would like to acknowledge the financial support of the School of Civil and Environmental Engineering, The University of New South Wales for his research fellowship for the period September 2012 onwards. 

\section*{References}
\bibliographystyle{elsarticle-num}
\bibliography{highorder}

\end{document}